\documentclass[10.5pt]{article}

\def\sqr#1#2{{\vcenter{\vbox{\hrule height.#2pt
              \hbox{\vrule width.#2pt height#1pt \kern#1pt \vrule width.#2pt}
              \hrule height.#2pt}}}}
\def\signed #1{{\unskip\nobreak\hfil\penalty50
              \hskip2em\hbox{}\nobreak\hfil#1
              \parfillskip=0pt \finalhyphendemerits=0 \par}}
\def\endpf{\signed {$\sqr69$}}

\def\dbR{{\mathop{\rm l\negthinspace R}}}

\def\3n{\negthinspace \negthinspace \negthinspace }
\def\2n{\negthinspace \negthinspace }
\def\1n{\negthinspace }

\def\dbH{{\mathop{\rm l\negthinspace H}}}

\def\dbR{{\mathop{\rm l\negthinspace R}}}
\def\={\buildrel \triangle \over =}

\def\ds{\displaystyle}

\def\ns{\noalign{\ss}}

\def\a{\alpha}
\def\b{\beta}
\def\g{\gamma}
\def\d{\delta}
\def\e{\varepsilon}

\def\k{\kappa}
\def\l{\lambda}
\def\m{\mu}
\def\n{\nu}
\def\si{\sigma}
\def\t{\tau}
\def\f{\varphi}
\def\th{\theta}
\def\o{\omega}

\def\Th{\Theta}

\def\cA{{\cal A}}

\def\cH{{\cal H}}

\def\cS{{\cal S}}

\def\cU{{\cal U}}

\def\no{\noindent}

\def\ss{\smallskip}
\def\ms{\medskip}
\def\bs{\bigskip}
\def\q{\quad}
\def\qq{\qquad}
\def\hb{\hbox}

\def\lan{\mathop{\langle}}
\def\ran{\mathop{\rangle}}

\def\rq{\eqno}
\def\pa{\partial}
\def\h{\widehat}
\def\wt{\widetilde}
\def\cd{\cdot}
\def\cds{\cdots}

\def\ae{\hbox{\rm a.e.{ }}}

\def\({\Big (}
\def\){\Big )}
\def\[{\Big[}
\def\]{\Big]}
\def\bde{\begin{definition}}
\def\ede{\end{definition}}
\def\be{\begin{equation}}
\def\bel{\begin{equation}\label}
\def\ee{\end{equation}}
\def\bt{\begin{theorem}}
\def\et{\end{theorem}}
\def\bc{\begin{corollary}}
\def\ec{\end{corollary}}
\def\bl{\begin{lemma}}
\def\el{\end{lemma}}
\def\bp{\begin{proposition}}
\def\ep{\end{proposition}}
\def\bas{\begin{assumption}}
\def\eas{\end{assumption}}
\def\br{\begin{remark}}
\def\er{\end{remark}}
\def\ba{\begin{array}}
\def\ea{\end{array}}
\def\ed{\end{document}}

\def\square#1{\vbox{\hrule\hbox{\vrule height#1%
     \kern#1\vrule}\hrule}}
\def\rectangle#1#2{\vbox{\hrule\hbox{\vrule height#1%
     \kern#2\vrule}\hrule}}

\font\tenbb=msbm10 \font\sevenbb=msbm7 \font\fivebb=msbm5

\newfam\bbfam
\scriptscriptfont\bbfam=\fivebb \textfont\bbfam=\tenbb
\scriptfont\bbfam=\sevenbb

\newtheorem{lemma}{Lemma}[section]
\newtheorem{remark}{Remark}[section]

\newtheorem{theorem}{Theorem}[section]
\newtheorem{corollary}{Corollary}[section]

\newtheorem{definition}{Definition}[section]
\newtheorem{proposition}{Proposition}[section]
\newtheorem{assumption}{Assumption}[section]

\makeatletter
   
   \@addtoreset{equation}{section}
\makeatother

\begin{document}

\title{\bf Hamilton-Jacobi Equations and Two-Person Zero-Sum
Differential Games with Unbounded Controls\footnote{This work is
supported in part by the NSF grant DMS-1007514, the NSFC grant
11171081, and the Postgraduate Scholarship Program of China.}}

\author{Hong Qiu$^{a,b}$~~and~~Jiongmin Yong$^{b}$
\\
 {\footnotesize\textsl{$^a$Department of Mathematics, Harbin
Institute of Technology, Weihai 264209, Shandong, China}} \\
{\footnotesize\textsl{$^b$Department of Mathematics, University of
Central Florida, Orlando, FL 32816, USA}}}

\maketitle

\begin{abstract}

A two-person zero-sum differential game with unbounded controls is
considered. Under proper coercivity conditions, the upper and lower
value functions are characterized as the unique viscosity solutions
to the corresponding upper and lower Hamilton--Jacobi--Isaacs
equations, respectively. Consequently, when the Isaacs' condition is
satisfied, the upper and lower value functions coincide, leading to
the existence of the value function of the differential game. Due to
the unboundedness of the controls, the corresponding upper and lower
Hamiltonians grow super linearly in the gradient of the upper and
lower value functions, respectively. A uniqueness theorem of
viscosity solution to Hamilton--Jacobi equations involving such kind
of Hamiltonian is proved, without relying on the convexity/concavity
of the Hamiltonian. Also, it is shown that the assumed coercivity
conditions guaranteeing the finiteness of the upper and lower value
functions are sharp in some sense.

\end{abstract}

\ms


\ms


\ms

\section{Introduction}

\ms

\rm

Let us begin with the following control system:

\bel{1.1}\left\{\ba{ll}
\ns\ds\dot y(s)=f(s,y(s),u_1(s),u_2(s)),\qq s\in[t,T],\\
\ns\ds y(t)=x.\ea\right.\ee
where $f:[0,T]\times\dbR^n\times U_1\times U_2\to\dbR^n$ is a given
map. In the above, $y(\cd)$ is the state trajectory taking values in
$\dbR^n$, and $(u_1(\cd),u_2(\cd))$ is the control pair taken from
the set $\cU_1^{\si_1}[t,T]\times\cU_2^{\si_2}[t,T]$ of {\it
admissible controls}, defined by the following:
$$\cU_i^{\si_i}[t,T]=\Big\{u_i:[t,T]\to U_i\Bigm|\|u_i(\cd)\|_{L^{\si_i}(t,T)}
\equiv\[\int_t^T|u_i(s)|^{\si_i}ds\]^{1\over\si_i}<\infty\Big\},\qq
i=1,2,$$
with $U_i$ being a closed subset of $\dbR^{m_i}$ and with some
$\si_i\ge1$. We point out that $U_1$ and $U_2$ are allowed to be
unbounded, and they could even be $\dbR^{m_1}$ and $\dbR^{m_2}$,
respectively. Hereafter, we suppress $\dbR^{m_i}$ in
$\|u_i(\cd)\|_{L^{\si_i}(t,T;\dbR^{m_i})}$ for notational simplicity
and this will not cause confusion. The {\it performance functional}
associated with (\ref{1.1}) is the following:
\bel{J}J(t,x;u_1(\cd),u_2(\cd))=\ds\int_t^Tg(s,y(s),u_1(s),u_2(s))ds+h(y(T)),\ee
with $g:[0,T]\times\dbR^n\times U_1\times U_2\to\dbR$ and
$h:\dbR^n\to\dbR$ being some given maps.

\ms

The above setting can be used to describe a two-person zero-sum
differential game: Player 1 wants to select a control
$u_1(\cd)\in\cU_1^{\si_1}[t,T]$ so that the functional (\ref{J}) is
minimized and Player 2 wants to select a control
$u_2(\cd)\in\cU_2^{\si_2}[t,T]$ so that the functional (\ref{J}) is
maximized. Therefore, $J(t,x;u_1(\cd),u_2(\cd))$ is a {\it cost
functional} for Player 1 and a {\it payoff functional} for Player 2,
respectively. If $U_2$ is a singleton, the above is reduced to a
standard optimal control problem.

\ms

Under some mild conditions, for any {\it initial pair}
$(t,x)\in[0,T]\times\dbR^n$ and control pair
$(u_1(\cd),u_2(\cd))\in\cU_1^{\si_1}[t,T]\times\cU_2^{\si_2}[t,T]$,
the state equation (\ref{1.1}) admits a unique solution
$y(\cd)\equiv y(\cd\,;t,x,u_1(\cd),u_2(\cd))$, and the performance
functional $J(t,x;u_1(\cd),u_2(\cd))$ is well-defined. By adopting
the notion of {\it Elliott--Kalton strategies} (\cite{Elliott-Kalton
1972}), we can define the {\it upper} and {\it lower value
functions} $V^{\pm}:[0,T]\times\dbR^n\to\dbR$ (see Section 3 for
details). Further, when $V^{\pm}(\cd\,,\cd)$ are differentiable,
they should satisfy the following upper and lower
Hamilton-Jacobi-Isaacs (HJI, for short) equations, respectively:
\bel{HJI0}\left\{\ba{ll}
\ns\ds
V^\pm_t(t,x)+H^\pm(t,x,V^\pm_x(t,x))=0,\qq(t,x)\in[0,T]\times\dbR^n,\\
\ns\ds V^\pm(T,x)=h(x),\qq x\in\dbR^n,\ea\right.\ee
where $H^\pm(t,x,p)$ are the so-called {\it upper} and {\it lower
Hamiltonians} defined by the following, respectively:
\bel{H}\left\{\ba{ll}
\ns\ds H^+(t,x,p)=\inf_{u_1\in U_1}\sup_{u_2\in
U_2}\[\lan p,f(t,x,u_1,u_2)\ran+g(t,x,u_1,u_2)\],\\
\ns\ds H^-(t,x,p)=\sup_{u_2\in U_2}\inf_{u_1\in U_1}\[\lan
p,f(t,x,u_1,u_2)\ran+g(t,x,u_1,u_2)\],\ea\right.
\q(t,x,p)\in[0,T]\times\dbR^n\times\dbR^n.\ee
When the sets $U_1$ and $U_2$ are bounded, the above differential
game is well-understood (\cite{Evans-Souganidis 1984,Ishii 1988}):
Under reasonable conditions, the upper and lower value functions
$V^\pm(\cd\,,\cd)$ are the unique viscosity solutions to the
corresponding upper and lower HJI equations, respectively.
Consequently, in the case that the following {\it Isaacs condition}:
\bel{Isaacs}H^+(t,x,p)=H^-(t,x,p),\qq\forall(t,x,p)\in[0,T]\times\dbR^n\times\dbR^n,\ee
holds, the upper and lower value functions coincide and the
two-person zero-sum differential game admits the value function
\bel{}V(t,x)=V^+(t,x)=V^-(t,x),\qq(t,x)\in[0,T]\times\dbR^n.\ee

\ms

For comparison purposes, let us now take a closer look at the
properties that the upper and lower value functions
$V^\pm(\cd\,,\cd)$ and the upper and lower Hamiltonians
$H^\pm(\cd\,,\cd\,,\cd)$ have, under classical assumptions. To this
end, let us recall the following classical assumption:

\ms

{\bf(B)} Functions $f:[0,T]\times\dbR^n\times U_1\times
U_2\to\dbR^n$, $g:[0,T]\times\dbR^n\times U_1\times U_2\to\dbR$, and
$h:\dbR^n\to\dbR$ are continuous. There exists a constant $L>0$ and
a continuous function $\o:[0,\infty)\times[0,\infty)\to[0,\infty)$,
increasing in each of its arguments and $\o(r,0)=0$ for all $r\ge0$,
such that for all $t,s\in[0,T],~x,y\in\dbR^n,~(u_1,u_2)\in U_1\times
U_2$,
\bel{1.6}\left\{\ba{ll}
\ns\ds|f(t,x,u_1,u_2)-f(s,y,u_1,u_2)|\le
L|x-y|+\o\big(|x|\vee|y|,|t-s|\big),\\
\ns\ds|g(t,x,u_1,u_2)-g(s,y,u_1,u_2)|\le\o\big(|x|\vee|y|,|x-y|+|t-s|\big),\\
\ns\ds|h(x)-h(y)|\le\o\big(|x|\vee|y|,|x-y|\big),\\
\ns\ds|f(t,0,u_1,u_2)|+|g(t,0,u_1,u_2)|+|h(0)|\le L,\ea\right.\ee
where $|x|\vee|y|=\max\{|x|,|y|\}$.

\ms

Condition (\ref{1.6}) implies that the continuity and the growth of
$(t,x)\mapsto(f(t,x,u_1,u_2)$, $g(t,x,u_1,u_2))$ are uniform in
$(u_1,u_2)\in U_1\times U_2$. This essentially will be the case if
$U_1$ and $U_2$ are bounded (or compact metric spaces). Let us state
the following proposition.

\ms

\bf Proposition 1.1. \sl Under assumption {\rm(B)}, one has the
following:

\ms

{\rm(i)} The upper and lower value functions $V^\pm(\cd\,,\cd)$ are
well-defined continuous functions. Moreover, they are the unique
viscosity solutions to the upper and lower HJI equations
$(\ref{HJI0})$, respectively. In particular, if Isaacs' condition
$(\ref{Isaacs})$ holds, the upper and lower value functions
coincide.

\ms

{\rm(ii)} The upper and lower Hamiltonians $H^\pm(\cd\,,\cd\,,\cd)$
satisfy the following: For all $t\in[0,T],~x,y,p,q\in\dbR^n$,
\bel{H-H}\left\{\ba{ll}
\ns\ds|H^\pm(t,x,p)-H^\pm(t,y,q)|\le
L(1+|x|)|p-q|+\o\big(|x|\vee|y|,|x-y|\big),\\
\ns\ds|H^\pm(t,x,p)|\le L(1+|x|)|p|+L+\o(|x|,|x|).\ea\right.\ee

\ms

\rm

Condition (\ref{H-H}) plays an important role in the proof of the
uniqueness of viscosity solution to HJI equations (\ref{HJI0})
(\cite{Crandall-Lions 1983,Ishii 1984}). Note that, in particular,
(\ref{H-H}) implies that $p\mapsto H^\pm(t,x,p)$ is at most of
linear growth.

\ms

Unfortunately, the above property (\ref{H-H}) fails, in general,
when the control domains $U_1$ and/or $U_2$ is unbounded. To make
this more convincing, let us look at a one-dimensional
linear-quadratic (LQ, for short) optimal control problem (which
amounts to saying that $U_1=\dbR$ and $U_2=\{0\}$). Consider the
state equation
$$\dot y(s)=y(s)+u(s),\qq s\in[t,T],$$
with a quadratic cost functional
$$J(t,x;u(\cd))={1\over2}\[\int_t^T\big(|y(s)|^2+|u(s)|^2\big)ds+|y(T)|^2\].$$
Then the Hamiltonian is
$$H(t,x,p)=\inf_{u\in\dbR}\[p(x+u)+{|x|^2+|u|^2\over2}\]=xp+{x^2\over2}-{p^2\over2}.$$
Thus, $p\mapsto H(t,x,p)$ is of quadratic growth and (\ref{H-H})
fails.

\ms

Optimal control problems with unbounded control domains were studied
in \cite{Bardi-Da Lio 1997, Da Lio 2000}. Uniqueness of viscosity
solution to the corresponding Hamilton-Jacobi-Bellman equation was
proved by some arguments relying on the convexity/concavity of the
corresponding Hamiltonian with respect to $p$. Recently, the above
results were substantially extended to stochastic optimal control
problems (\cite{Da Lio--Ley 2011}). On the other hand, as an
extension of \cite{Yong 1994}, two-person zero-sum differential
games with (only) one player having unbounded control were studied
in \cite{Rampazzo 1998}. Some nonlinear $H_\infty$ problems can also
be treated as such kind of differential games \cite{McEneaney
1998,Soravia 1999}. Further, stochastic two-person zero-sum
differential games were studied in \cite{Da Lio--Ley 2006} with one
player having unbounded control and with the two players' controls
being separated both in the state equation and the performance
functional.

\ms

The main purpose of this paper is to study two-person zero-sum
differential games with both players having unbounded controls, and
the controls of two players are not necessarily separated. One
motivation comes from the problem of what we call the
affine-quadratic (AQ, for short) two-person zero-sum differential
games, by which we mean that the right hand side of the state
equation is affine in the controls, and the integrand of the
performance functional is quadratic in the controls (see Section 2).
This is a natural generalization of the classical LQ problems. For
general two-person zero-sum differential games with (both players
having) unbounded controls, under some mild coercivity conditions,
the upper and lower Hamiltonians $H^\pm(t,x,p)$ are proved to be
well-defined, continuous, and locally Lipschitz in $p$. Therefore,
the upper and lower HJI equations can be formulated. Then we will
establish the uniqueness of viscosity solutions to a general first
order Hamilton-Jacobi equation which includes our upper and lower
HJI equations of the differential game. Comparing with a relevant
result found in \cite{Crandall-Lions 1986}, the conditions we
assumed here are a little different from theirs and we present a
detailed proof for reader's convenience. By assuming a little
stronger coercivity conditions, together with some additional
conditions (guaranteeing the well-posedness of the state equation,
etc.), we show that the upper and lower value functions can be
well-defined and are continuous. Combining the above results, one
obtains a characterization of the upper and lower value functions of
the differential game as the unique viscosity solutions to the
corresponding upper and lower HJI equations. Then if in addition,
the Isaacs' condition holds, the upper and lower value functions
coincide which yields the existence of the value function of the
differential game.

\ms

We would like to mention here that due to the unboundedness of the
controls, the continuity of the upper and lower value functions
$V^\pm(t,x)$ in $t$ is quite subtle. To prove that, we need to
establish a modified principle of optimality and fully use the
coercivity conditions. It is interesting to indicate that the
assumed coercivity conditions that ensuring the finiteness of the
upper and lower value functions are actually sharp in some sense,
which was illustrated by a one-dimensional LQ situation.

\ms

For some other relevant works in the literature, we would like to
mention \cite{Lions-Souganidis 1985,Friedman-Soudanidis
1986,Fleming-Souganidis 1989,Bardi-Capuzzo-Dolcetta 1997,You 2002,
Soravia 2004}, and references cited therein.

\ms

The rest of the paper is organized as follows. In Section 2, we make
some brief observations on an AQ two-person differential game, for
which we have a situation that the Isaacs' condition holds and the
upper and lower Hamiltonians $H^\pm(t,x,p)$ are quadratic in $p$ but
may be neither convex nor concave. Section 3 is devoted to a study
of upper and lower Hamiltonians. The uniqueness of viscosity
solutions to a class of HJ equations will be proved in Section 4. In
Section 5, we will show that under certain conditions, the upper and
lower value functions are well-defined and continuous. Finally, in
Section 6, we show that the assumed coercivity conditions ensuring
the upper and lower value functions to be well-defined are sharp in
some sense.

\ms

\section{An Affine-Quadratic Two-Person Differential Game}

To better understand two-person zero-sum differential games with
unbounded controls, in this section, we look at a nontrivial special
case which is a main motivation of this paper. Consider the
following state equation:
\bel{2.1}\left\{\ba{ll}
\ns\ds\dot y(s)=A(s,y(s))+B_1(s,y(s))u_1(s)+B_2(s,y(s))u_2(s),\qq
s\in[t,T],\\
\ns\ds y(t)=x,\ea\right.\ee
for some suitable matrix valued functions $A(\cd\,,\cd)$,
$B_1(\cd\,,\cd)$, and $B_2(\cd\,,\cd)$. The state $y(\cd)$ takes
values in $\dbR^n$ and the control $u_i(\cd)$ takes values in
$U_i=\dbR^{m_i}$ ($i=1,2$). The performance functional is given by
\bel{}\ba{ll}
\ns\ds J(t,x;u_1(\cd),u_2(\cd))=\int_t^T\[Q(s,y(s))+{1\over2}\lan
R_1(s,y(s))u_1(s),u_1(s)\ran\\
\ns\ds\qq\qq\qq\qq\qq\q+\lan
S(s,y(s))u_1(s),u_2(s)\ran-{1\over2}\lan
R_2(s,y(s))u_2(s),u_2(s)\ran\\
\ns\ds\qq\qq\qq\qq\qq\q+\lan\th_1(s,y(s)),u_1(s)\ran+\lan\th_2(s,y(s)),u_2(s)\ran\]ds+G(y(T)),\ea\ee
for some scalar functions $Q(\cd\,,\cd)$ and $G(\cd)$, some vector
valued functions $\th_1(\cd\,,\cd)$ and $\th_2(\cd\,,\cd)$, and some
matrix valued functions $R_1(\cd\,,\cd)$, $R_2(\cd\,,\cd)$, and
$S(\cd\,,\cd)$. Note that the right hand side of the state equation
is affine in the controls $u_1(\cd)$ and $u_2(\cd)$, and the
integrand in the performance functional is up to quadratic in
$u_1(\cd)$ and $u_2(\cd)$. Therefore, we refer to such a problem as
an {\it affine-quadratic} (AQ, for short) {\it two-person zero-sum
differential game}. We also note that due to the presence of the
term $\lan S(s,y(s))u_1(s),u_2(s)\ran$, controls $u_1(\cd)$ and
$u_2(\cd)$ cannot be completely separated. Let us now introduce the
following basic hypotheses concerning the above AQ two-person
zero-sum differential game.

\ms

{\bf(AQ1)} The maps
$$\ba{ll}
\ns\ds A:[0,T]\times\dbR^n\to\dbR^n,\q
B_1:[0,T]\times\dbR^n\to\dbR^{n\times m_1},\q
 B_2:[0,T]\times\dbR^n\to\dbR^{n\times m_2},\ea$$
are continuous.

\ms

{\bf(AQ2)} The maps
$$\ba{ll}
\ns\ds Q:[0,T]\times\dbR^n\to\dbR,\q G:\dbR^n\to\dbR,\q
R_1:[0,T]\times\dbR^n\to\cS^{m_1},\q
R_2:[0,T]\times\dbR^n\to\cS^{m_2},\\
\ns\ds S:[0,T]\times\dbR^n\to\dbR^{m_2\times m_1},\q
\th_1:[0,T]\times\dbR^n\to\dbR^{m_1},\q\th_2:[0,T]\times\dbR^n\to\dbR^{m_2}\ea$$
are continuous (where $\cS^m$ stands for the set of all $(m\times
m)$ symmetric matrices), and $R_1(t,x)$ and $R_2(t,x)$ are positive
definite for all $(t,x)\in[0,T]\times\dbR^n$.

\ms

With the above hypotheses, we let
\bel{}\ba{ll}
\ns\ds\dbH(t,x,p,u_1,u_2)=\lan
p,A(t,x)+B_1(t,x)u_1+B_2(t,x)u_2\ran+Q(t,x)+{1\over2}\lan
R_1(t,x)u_1,u_1\ran\\
\ns\ds\qq\qq\qq\qq+\lan S(t,x)u_1,u_2\ran-{1\over2}\lan
R_2(t,x)u_2,u_2\ran+\lan\th_1(t,x),u_1\ran+\lan\th_2(t,x),u_2\ran.\ea\ee
Our result concerning the above-defined function is the following
proposition.

\ms

\bf Proposition 2.1. \rm Let {\rm(AQ1)--(AQ2)} hold. Then the matrix
$\pmatrix{R_1(t,x)&S(t,x)^T\cr S(t,x)&-R_2(t,x)}$ is invertible, and
\bel{dbH}\ba{ll}
\ns\ds\dbH(t,x,p,u_1,u_2)={1\over2}\lan R_1(t,x)(u_1-\bar
u_1),u_1-\bar
u_1\ran+\lan S(t,x)(u_1-\bar u_1),u_2-\bar u_2\ran\\
\ns\ds\qq\qq\qq\qq\q-{1\over2}\lan R_2(t,x)(u_2-\bar u_2),u_2-\bar
u_2\ran+Q_0(t,x,p),\ea\ee
where
\bel{u}\pmatrix{\bar u_1\cr\bar u_2}=-\pmatrix{R_1(t,x)&S(t,x)^T\cr
S(t,x)&-R_2(t,x)}^{-1}\pmatrix{B_1(t,x)^Tp+\th_1(t,x)\cr
B_2(t,x)^Tp+\th_2(t,x)},\ee
and
\bel{Q0}\ba{ll}
\ns\ds
Q_0(t,x,p)=Q(t,x)+\lan p,A(t,x)\ran\\
\ns\ds\qq\qq\qq-{1\over2}\pmatrix{B_1(t,x)^Tp+\th_1(t,x)\cr
B_2(t,x)^Tp+\th_2(t,x)}^T\2n\pmatrix{R_1(t,x)&S(t,x)^T\cr
S(t,x)&-R_2(t,x)}^{-1}\2n\pmatrix{B_1(t,x)^Tp+\th_1(t,x)\cr
B_2(t,x)^Tp+\th_2(t,x)},\ea\ee
Further, $(\bar u_1,\bar u_2)$ given by (\ref{u}) is the unique
saddle point of $(u_1,u_2)\mapsto\dbH(t,x,p,u_1,u_2)$, namely,
\bel{}\dbH(t,x,p,\bar u_1,u_2)\le\dbH(t,x,p,\bar u_1,\bar
u_2)\le\dbH(t,x,p,u_1,\bar u_2),\qq\forall(u_1,u_2)\in U_1\times
U_2,\ee
and consequently, the Isaacs' condition is satisfied:
\bel{H=H}\ba{ll}
\ns\ds H^+(t,x,p)\equiv\inf_{u_1\in U_1}\sup_{u_2\in
U_2}\dbH(t,x,p,u_1,u_2)=\sup_{u_2\in U_2}\inf_{u_1\in
U_1}\dbH(t,x,p,u_1,u_2)\\
\ns\ds\qq\qq\q\equiv
H^-(t,x,p)=Q_0(t,x,p),\qq\qq\forall(t,x,p)\in[0,T]\times\dbR^n\times\dbR^n.\ea
\ee

\ms

\it Proof. \rm For simplicity of notation, let us suppress $(t,x)$
below. We may write
$$\ba{ll}
\ns\ds\dbH(p,u_1,u_2)={1\over2}\lan R_1(u_1-\bar u_1),u_1-\bar
u_1\ran+\lan S(u_1-\bar u_1),u_2-\bar u_2\ran-{1\over2}\lan
R_2(u_2-\bar u_2),u_2-\bar u_2\ran+Q_0,\ea$$
with $\bar u_1$, $\bar u_2$, and $Q_0$ undetermined. Then
$$\ba{ll}
\ns\ds\lan p,A\ran+Q+\lan B_1^Tp+\th_1,u_1\ran+\lan
B_2^Tp+\th_2,u_2\ran+{1\over2}\lan R_1u_1,u_1\ran+\lan
Su_1,u_2\ran-{1\over2}\lan
R_2u_2,u_2\ran\\
\ns\ds=\dbH(p,u_1,u_2)={1\over2}\lan R_1u_1,u_1\ran+\lan
Su_1,u_2\ran-{1\over2}\lan R_2u_2,u_2\ran-\lan R_1\bar
u_1,u_1\ran\\
\ns\ds\q-\lan S^T\bar u_2,u_1\ran-\lan S\bar u_1,u_2\ran+\lan
R_2\bar u_2,u_2\ran+{1\over2}\lan R_1\bar u_1,\bar u_1\ran+\lan
S\bar u_1,\bar u_2\ran-{1\over2}\lan R_2\bar u_2,\bar
u_2\ran+Q_0.\ea$$
Hence, we must have
\bel{2.9}\left\{\ba{ll}
\ns\ds B_1^Tp+\th_1=-R_1\bar u_1-S^T\bar u_2,\qq B_2^Tp+\th_2=-S\bar u_1+R_2\bar u_2,\\
\ns\ds\lan p,A\ran+Q={1\over2}\lan R_1\bar u_1,\bar u_1\ran+\lan
S\bar u_1,\bar u_2\ran-{1\over2}\lan R_2\bar u_2,\bar
u_2\ran+Q_0.\ea\right.\ee
Consequently, from the first two equations in (\ref{2.9}), we have
$$\pmatrix{R_1&S^T\cr S&-R_2}\pmatrix{\bar u_1\cr\bar
u_2}=-\pmatrix{B_1^Tp+\th_1\cr B_2^Tp+\th_2}.$$
Note that
$$\ba{ll}
\ns\ds\det\pmatrix{R_1&S^T\cr
S&-R_2}=\det\pmatrix{R_1&0\cr0&-(R_2+SR_1^{-1}S^T)}=(-1)^{m_2}\det(R_1)\det(R_2+SR_1^{-1}S^T)\ne0.\ea$$
Thus, $\pmatrix{R_1&S^T\cr S&-R_2}$ is invertible, which yields
$$\pmatrix{\bar u_1\cr\bar u_2}=-\pmatrix{R_1&S^T\cr
S&-R_2}^{-1}\pmatrix{B_1^Tp+\th_1\cr B_2^Tp+\th_2}.$$
Then from the last equality in (\ref{2.9}), one has
$$\ba{ll}
\ns\ds Q_0=\lan p,A\ran+Q-{1\over2}\pmatrix{\bar u_1\cr\bar
u_2}^T\pmatrix{R_1&S^T\cr S&-R_2}\pmatrix{\bar u_1\cr\bar
u_2}\\
\ns\ds\q\;\,=\lan p, A\ran+Q-{1\over2}\pmatrix{B_1^Tp+\th_1\cr
B_2^Tp+\th_2}^T\pmatrix{R_1&S^T\cr
S&-R_2}^{-1}\pmatrix{B_1^Tp+\th_1\cr B_2^Tp+\th_2},\ea$$
proving (\ref{dbH}). Now, we see that
$$\ba{ll}
\ns\ds\dbH(p,\bar u_1,u_2)=-{1\over2}\lan R_2(u_2-\bar
u_2),u_2-\bar u_2\ran+Q_0\le Q_0=\dbH(p,\bar u_1,\bar u_2)\\
\ns\ds\qq\qq\q\,\le{1\over2}\lan R_1(u_1-\bar u_1),u_1-\bar
u_1\ran+Q_0(t,x,p)=\dbH(p,u_1,\bar u_2),\ea$$
which means that $(\bar u_1,\bar u_2)$ is a saddle point of
$\dbH(t,x,p,u_1,u_2)$. Then the Isaacs condition (\ref{H=H}) follows
easily. Finally, since $R_1$ and $R_2$ are positive definite, the
saddle point must be unique. \endpf

\ms

We see that in the current case, $p\mapsto H^\pm(t,x,p)$ is
quadratic, and is neither convex nor concave in general. As a matter
of fact, the Hessian $H^\pm_{pp}(t,x,p)$ of $H^\pm(t,x,p)$ is given
by the following:
$$H^\pm_{pp}(t,x,p)=-{1\over2}\pmatrix{B_1(t,x)^T\cr
B_2(t,x)^T}^T\2n\pmatrix{R_1(t,x)&S(t,x)^T\cr
S(t,x)&-R_2(t,x)}^{-1}\2n\pmatrix{B_1(t,x)^T\cr B_2(t,x)^T}.$$
which is indefinite in general.

\ms

We have seen from the above that in order the upper and lower
Hamiltonians to be well-defined, the only crucial assumption that we
made is the positive definiteness of the matrix-valued maps
$R_1(\cd\,,\cd)$ and $R_2(\cd\,,\cd)$. Whereas, in order to study
the AQ two-person zero-sum differential games, we need a little
stronger hypotheses. For example, in order the state equation to be
well-posed, we need the right hand side of the state equation is
Lipschitz continuous in the state variable, for any given pair of
controls, etc. We will look at the general situation a little later.

\section{Upper and Lower Hamiltonians}

In this section, we will carefully look at the upper and lower
Hamiltonians associated with general two-person zero-sum
differential games with unbounded controls. First of all, we
introduce the following standing assumptions.

\ms

{\bf(H0)} For $i=1,2$, the set $U_i\subseteq\dbR^{m_i}$ is closed
and
\bel{0}0\in U_i,\qq i=1,2.\ee
The time horizon $T>0$ is fixed.

\ms

Note that both $U_1$ and $U_2$ could be unbounded and may even be
equal to $\dbR^{m_1}$ and $\dbR^{m_2}$, respectively. Condition
(\ref{0}) is for convenience. We may make a translation of the
control domains and make corresponding changes in the control
systems and performance functional to achieve this.


\ms

Inspired by the AQ two-person zero-sum differential games, let us
now introduce the following assumptions for the involved functions
$f$ and $g$ in the state equation (\ref{1.1}) and the performance
functional (\ref{J}). We denote $\lan x\ran=\sqrt{1+|x|^2}$.

\ms

{\bf(H1)} Map $f:[0,T]\times\dbR^n\times U_1\times U_2\to\dbR^n$ is
continuous and there are constants $\si_1,\si_2\ge0$ such that
\bel{|f|}\ba{ll}
\ns\ds|f(t,x,u_1,u_2)|\le L\big(\1n\lan
x\ran+|u_1|^{\si_1}+|u_2|^{\si_2}\big),
\qq\forall(t,x,u_1,u_2)\in[0,T]\times\dbR^n\times U_1\times
U_2.\ea\ee

\ms

{\bf(H2)} Map $g:[0,T]\times\dbR^n\times U_1\times U_2\to\dbR$ is
continuous and there exist constants $L,c,\rho_1,\rho_2>0$ and
$\m\ge1$ such that
\bel{g}\ba{ll}
\ns\ds c|u_1|^{\rho_1}-L\big(\1n\lan
x\ran{}^\m+|u_2|^{\rho_2}\big)\le g(t,x,u_1,u_2)\le L\big(\1n\lan
x\ran{}^{\m}+|u_1|^{\rho_1}\big)-c|u_2|^{\rho_2},\\
\ns\ds\qq\qq\qq\qq\qq\qq\qq\forall(t,x,u_1,u_2)\in
[0,T]\times\dbR^n\times U_1\times U_2.\ea\ee

\ms

Further, we introduce the following compatibility condition which
will be crucial below.

\ms

{\bf(H3)} The constants $\si_1,\si_2,\rho_1,\rho_2$ in (H1)--(H2)
satisfy the following:
\bel{}\si_i<\rho_i,\qq i=1,2.\ee

It is not hard to see that the above (H1)--(H3) includes the AQ
two-person zero-sum differential game described in the previous
section as a special case. Now, we let
\bel{dbH}\ba{ll}
\ns\ds\dbH(t,x,p,u_1,u_2)=\lan
p,f(t,x,u_1,u_2)\ran+g(t,x,u_1,u_2),\qq(t,x,u_1,u_2)\in[0,T]\times
\dbR^n\times U_1\times U_2.\ea\ee
Then the {\it upper} and {\it lower Hamiltonians} are defined as
follows:
\bel{Hpm}\left\{\ba{ll}
\ns\ds H^+(t,x,p)=\inf_{u_1\in U_1}\sup_{u_2\in
U_2}\dbH(t,x,p,u_1,u_2),\\
\ns\ds H^-(t,x,p)=\sup_{u_2\in U_2}\inf_{u_1\in
U_1}\dbH(t,x,p,u_1,u_2),\ea\right.\qq(t,x,p)\in[0,T]\times\dbR^n\times\dbR^n,\ee
provided the involved infimum and supremum exist. Note that the
upper and lower Hamiltonians are nothing to do with the function
$h(\cd)$ (appears as the terminal cost/payoff in (\ref{J})). The
main result of this section is the following.

\ms

\bf Proposition 3.1. \sl Under {\rm(H1)--(H3)}, the upper and lower
Hamitonians $H^\pm(\cd\,,\cd\,,\cd)$ are well-defined and
continuous. Moreover, there are constants $C>0$, $\l_i,\n_i\ge0$,
$(i=1,2,\cds,k)$ such that
\bel{|H|}\ba{ll}
\ns\ds-L\lan x\ran{}^\m-L\lan
x\ran|p|-C|p|^{\rho_1\over\rho_1-\si_1} \le H^\pm(t,x,p)\le L\lan
x\ran{}^\m+L\lan x\ran|p|+C|p|^{\rho_2\over\rho_2-\si_2},\\
\ns\ds\qq\qq\qq\qq\qq\qq\qq\qq\qq\forall(t,x,p)\in[0,T]\times\dbR^n\times\dbR^n,\ea\ee
and
\bel{|H-H|}\ba{ll}
\ns\ds|H^\pm(t,x,p)-H^\pm(t,x,q)|\le C\sum_{i=1}^k\lan
x\ran{}^{\l_i}\big(|p|\vee|q|\big)^{\n_i}|p-q|,\\
\ns\ds\qq\qq\qq\qq\qq\qq\qq\qq\forall(t,x)\in[0,T]\times\dbR^n,~p,q\in\dbR^n.\ea\ee

\ms

\rm

To prove the above, we will use the following lemma.

\ms

\bf Lemma 3.2. \sl Let $0<\si<\rho$ and $c,N>0$. Let
$$\th(r)=Nr^\si-cr^\rho,\qq r\in[0,\infty).$$
Then
\bel{}\ba{ll}
\ns\ds\max_{r\in[0,\infty)}\th(r)=\max_{r\in[0,\bar
r]}\th(r)=\th(\bar r)=(\rho-\si)\({\si^\si N^\rho\over\rho^\rho
c^\si}\)^{1\over\rho-\si},\ea\ee
with
\bel{bar r}\bar r=\({\si N\over \rho c}\)^{1\over\rho-\si}\,.\ee

\it Proof. \rm From
$$\th(0)=0,\q\lim_{r\to\infty}\th(r)=-\infty,$$
we see that the maximum of $\th(\cd)$ on $[0,\infty)$ is achieved at
some point $\bar r\in(0,\infty)$. Set
$$0=\th'(r)=N\si r^{\si-1}-c\rho r^{\rho-1}.$$
Then
$$r^{\rho-\si}={N\si\over c\rho}>0\,,$$
which implies that the maximum is achieved at $\bar r$ given by
(\ref{bar r}), and
$$\ba{ll}
\ns\ds\max_{r\in[0,\infty)}\th(r)=\max_{r\in[0,\bar
r]}\th(r)=\th(\bar r)=N\({N\si\over
c\rho}\)^{\si\over\rho-\si}-c\({N\si\over
c\rho}\)^{\rho\over\rho-\si}\\
\ns\ds\qq\qq\q\,=\[\({\si\over
c\rho}\)^{\si\over\rho-\si}-c\({\si\over
c\rho}\)^{\rho\over\rho-\si}\]N^{\rho\over\rho-\si}
={(\rho-\si)\si^{\si\over\rho-\si}\over
c^{\si\over\rho-\si}\rho^{\rho\over\rho-\si}}N^{\rho\over\rho-\si}.\ea$$
This proves our conclusion. \endpf

\ms

\it Proof of Proposition 3.1. \rm Let us look at $H^+(t,x,p)$
carefully ($H^-(t,x,p)$ can be treated similarly). First, by our
assumption, we have
\bel{dbH<}\ba{ll}
\ns\ds\dbH(t,x,p,u_1,u_2)\le|p|\,|f(t,x,u_1,u_2)|+g(t,x,u_1,u_2)\\
\ns\ds\le L\big(\1n\lan
x\ran+|u_1|^{\si_1}+|u_2|^{\si_2}\big)|p|+L\big(\1n\lan
x\ran{}^\m+|u_1|^{\rho_1}\big)-c|u_2|^{\rho_2}\\
\ns\ds=L\big(\1n\lan x\ran{}^\m+\lan
x\ran|p|+|p|\,|u_1|^{\si_1}+|u_1|^{\rho_1}\big)+L|p|\,|u_2|^{\si_2}-c|u_2|^{\rho_2},\ea\ee
and
\bel{dbH>}\ba{ll}
\ns\ds\dbH(t,x,p,u_1,u_2)\ge-|p|\,|f(t,x,u_1,u_2)|+g(t,x,u_1,u_2)\\
\ns\ds\ge-L\big(\1n\lan
x\ran{}+|u_1|^{\si_1}+|u_2|^{\si_2}\big)|p|-L\big(\1n\lan
x\ran{}^\m+|u_2|^{\rho_2}\big)+c|u_1|^{\rho_2}\\
\ns\ds=\1n-L\big(\1n\lan x\ran{}^\m\1n+\1n\lan
x\ran|p|\1n+\1n|p|\,|u_2|^{\si_2}\1n+\1n|u_2|^{\rho_2}\big)\1n-\1n
L|p|\,|u_1|^{\si_1}\1n+\1n c|u_1|^{\rho_1}.\ea\ee
Noting $\si_1<\rho_1$, from (\ref{dbH<}), we see that for any fixed
$(t,x,p,u_1)\in[0,T]\times\dbR^n\times\dbR^n\times U_1$, the map
$u_2\mapsto\dbH(t,x,p,u_1,u_2)$ is coercive from above.
Consequently, since $U_2$ is closed, for any given
$(t,x,p,u_1)\in[0,T]\times\dbR^n\times\dbR^n\times U_1$, there
exists a $\bar u_2\equiv\bar u_2(t,x,p,u_1)\in U_2$ such that
\bel{cH+}\ba{ll}
\ns\ds\cH^+(t,x,p,u_1)\equiv\sup_{u_2\in
U_2}\dbH(t,x,p,u_1,u_2)=\sup_{u_2\in U_2,\,|u_2|\le|\bar
u_2|}\dbH(t,x,p,u_1,u_2)=\dbH(t,x,p,u_1,\bar u_2)\\
\ns\ds\le L\big(\1n\lan x\ran{}^\m+\lan x\ran|p|+|p|\,
|u_1|^{\si_1}+|u_1|^{\rho_1}\big)+L|p|\,|\bar u_2|^{\si_2}
-c|\bar u_2|^{\rho_2}\\
\ns\ds\le\1n L\big(\1n\lan x\ran{}^\m\2n+\1n\1n\lan
x\ran|p|\1n+\1n|p|\,
|u_1|^{\si_1}\1n+\1n|u_1|^{\rho_1}\big)+(\rho_2-\si_2)\({\si_2^{\si_2}
\big(\1n L|p|\big)^{\rho_2}
\over\rho_2^{\rho_2}c^{\si_2}}\)^{1\over\rho_2-\si_2}\\
\ns\ds\le L\big(\1n\lan x\ran{}^\m+\lan
x\ran|p|+|p|\,|u_1|^{\si_1}+|u_1|^{\rho_1}\big)
+K_2|p|^{\rho_2\over\rho_2-\si_2},\ea\ee
where
$$K_2=(\rho_2-\si_2)\({\si_2^{\si_2}L^{\rho_2}\over\rho_2^{\rho_2}c^{\si_2}}\)^{1\over\rho_2-\si_2}.$$
Here, we have used Lemma 3.2. On the other hand, from (\ref{dbH>}),
for any $(t,x,p,u_1)\in[0,T]\times\dbR^n\times\dbR^n\times U_1$, we
have
\bel{sup H+>}\ba{ll}
\ns\ds\cH^+(t,x,p,u_1)=\sup_{u_2\in
U_2}\dbH(t,x,p,u_1,u_2)\ge\dbH(t,x,p,u_1,0)\\
\ns\ds\qq\qq\qq\;\ge-L\big(\1n\lan x\ran{}^\m+\lan
x\ran|p|\big)-L|p|\,|u_1|^{\si_1}+c|u_1|^{\rho_1}.\ea\ee
By Young's inequality, we have
$$L|p|\,|u_i|^{\si_i}\le{c\over2}|u_i|^{\rho_i}+\bar K_i|p|^{\rho_i\over\rho_i-\si_i},\qq
i=1,2,$$
for some absolute constants $\bar K_i$ (depending on
$L,c,\rho_i,\si_i$ only), which leads to
\bel{3.16}{c\over2}|u_i|^{\rho_i}\le
c|u_i|^{\rho_i}-L|p|\,|u_i|^{\si_i}+\bar
K_i|p|^{\rho_i\over\rho_i-\si_i},\qq i=1,2.\ee
Hence, combining the first inequality in (\ref{cH+}) and (\ref{sup
H+>}), we obtain
\bel{bar u2}\ba{ll}
\ns\ds{c\over2}|\bar u_2|^{\rho_2}\le c|\bar
u_2|^{\rho_2}-L|p|\,|\bar u_2|^{\si_2}+\bar
K_2|p|^{\rho_2\over\rho_2-\si_2}\\
\ns\ds\le\1n L\big(\1n\lan x\ran{}^\m\1n+\1n\lan
x\ran|p|\1n+\1n|p|\,|u_1|^{\si_1}\1n+\1n|u_1|^{\rho_1}\big)\1n
-\1n\cH^+(t,x,u_1)\1n+\1n\bar K_2|p|^{\rho_2\over\rho_2-\si_2}\\
\ns\ds\le2L\big(\1n\lan x\ran{}^\m+\lan
x\ran|p|+|p|\,|u_1|^{\si_1}\big)+(L-c)|u_1|^{\rho_1} +\bar
K_2|p|^{\rho_2\over\rho_2-\si_2}\equiv\h
K_2\big(|x|,|p|,|u_1|\big).\ea\ee
The above implies that for any compact set
$G\subseteq[0,T]\times\dbR^n\times\dbR^n\times U_1$, there exists a
compact set $\h U_2(G)\subseteq U_2$, depending on $G$, such that
$$\cH^+(t,x,p,u_1)=\sup_{u_2\in\h
U_2(G)}\dbH(t,x,p,u_1,u_2),\qq\forall(t,x,p,u_1)\in G.$$
Hence, $\cH^+(\cd\,,\cd\,,\cd\,,\cd)$ is continuous. Next, from
(\ref{sup H+>}), noting $\si_1<\rho_1$, we have that for any fixed
$(t,x,p)\in[0,T]\times\dbR^n\times\dbR^n$, the map
$u_1\mapsto\cH^+(t,x,p,u_1)$ is coercive from below. Therefore,
using the continuity of $\cH^+(\cd\,,\cd\,,\cd\,,\cd)$, one can find
a $\bar u_1\equiv\bar u_1(t,x,p)$ such that (note (\ref{3.16}))
\bel{sup H>}\ba{ll}
\ns\ds H^+(t,x,p)=\inf_{u_1\in U_1}\sup_{u_2\in
U_2}\dbH(t,x,p,u_1,u_2)=\inf_{u_1\in U_1}\cH^+(t,x,p,u_1)=\cH^+(t,x,p,\bar u_1)\\
\ns\ds\ge\inf_{u_1\in U_1}\dbH(t,x,p,u_1,0)\ge\inf_{u_1\in
U_1}\Big\{-L\big(\1n\lan x\ran{}^\m+\lan x\ran|p|\big)
-L|p|\,|u_1|^{\si_1}+c|u_1|^{\rho_1}\Big\}\\
\ns\ds\ge-L\big(\1n\lan x\ran{}^\m+\lan x\ran|p|\big)-\bar
K_1|p|^{\rho_1\over\rho_1-\si_1}.\ea\ee
This means that $H^+(t,x,p)$ is well-defined for all
$(t,x,p)\in[0,T]\times\dbR^n\times\dbR^n$, and it is locally bounded
from below. Also, from (\ref{cH+}), we obtain
\bel{H+<}\ba{ll}
\ns\ds H^+(t,x,p)=\inf_{u_1\in U_1}\sup_{u_2\in
U_2}\dbH(t,x,p,u_1,u_2)\le\sup_{u_2\in
U_2}\dbH(t,x,p,0,u_2)\equiv\cH^+(t,x,p,0)\\
\ns\ds\qq\qq\q\le L\big(\1n\lan x\ran{}^\m+\lan
x\ran|p|\big)+K_2|p|^{\rho_2\over\rho_2-\si_2}. \ea\ee
This proves (\ref{|H|}) for $H^+(\cd\,,\cd\,,\cd)$.

\ms

Next, we want to get the local Lipschitz continuity of the map
$p\mapsto H^+(t,x,p)$. To this end, we first let
$$\ba{ll}
\ns\ds U_1(|x|,|p|)=\Big\{u_1\in
U_1\bigm|{c\over2}|u_1|^{\rho_1}\le2L\big(\1n\lan x\ran{}^\m+\lan
x\ran|p|\big)+\bar K_1|p|^{\rho_1\over\rho_1-\si_1}+K_2
|p|^{\rho_2\over\rho_2-\si_2}+1\Big\},\q\forall x,p\in\dbR^n,\ea$$
which, for any given $x,p\in\dbR^n$, is a compact set. Clearly, for
any $u_1\in U_1\setminus U_1(|x|,|p|)$, one has (note (\ref{3.16}))
$$c|u_1|^{\rho_1}-L|p|\,|u_1|^{\si_1}\ge{c\over2}|u_1|^{\rho_1}-\bar
K_1|p|^{\rho_1\over\rho_1-\si_1}>2L\big(\1n\lan x\ran{}^\m+\lan
x\ran|p|\big)+K_2|p|^{\rho_2\over\rho_2-\si_2}+1.$$
Thus, for such a $u_1$, by (\ref{sup H+>}) and (\ref{H+<}),
\bel{}\ba{ll}
\ns\ds\cH^+(t,x,p,u_1)\ge-L\big(\1n\lan x\ran{}^\m+\lan
x\ran|p|\big)-L|p|\,|u_1|^{\si_1}+c|u_1|^{\rho_1}\\
\ns\ds>L\big(\1n\lan x\ran{}^\m+\lan x\ran|p|\big)
+K_2|p|^{\rho_2\over\rho_2-\si_2}+1\ge H^+(t,x,p)+1=\inf_{u_1\in
U_1}\cH^+(t,x,p,u_1)+1.\ea\ee
Hence,
\bel{}\inf_{u_1\in U_1}\cH^+(t,x,p,u_1)=\inf_{u_1\in
U_1(|x|,|p|)}\cH^+(t,x,p,u_1).\ee
Now, for any $u_1\in U_1(|x|,|p|)$, by (\ref{bar u2}), we have
\bel{}\ba{ll}
\ns\ds{c\over2}|\bar u_2|^{\rho_2}\le\h
K_2\big(|x|,|p|,|u_1|\big)\le\wt K_2(|x|,|p|),\ea\ee
for some $\wt K_2(|x|,|p|)$. Hence, if we let
$$U_2(|x|,|p|)=\Big\{u_2\in U_2\bigm|{c\over2}|u_2|^{\rho_2}\le\wt
K_2(|x|,|p|)\Big\},$$
which is a compact set (for any given $x,p\in\dbR^n$), then for any
$(t,x,p)\in[0,T]\times\dbR^n\times\dbR^n$,
\bel{3.15}H^+(t,x,p)=\inf_{u_1\in U_1(|x|,|p|)}\sup_{u_2\in
U_2(|x|,|p|)}\dbH(x,p,u_1,u_2).\ee
This implies that $H^+(\cd\,,\cd\,,\cd)$ is continuous. Next, we
look at some estimates. By definition, for any $u_1\in
U_1(|x|,|p|)$, we have
$$\ba{ll}|u_1|^{\rho_1}\le C\Big\{\lan
x\ran{}^\m+\lan
x\ran|p|+|p|^{\rho_1\over\rho_1-\si_1}+|p|^{\rho_2\over\rho_2-\si_2}\Big\}.\ea$$
Therefore,
$$\ba{ll}
\ns\ds|u_1|^{\si_1}\le C\Big\{\lan x\ran{}^\m+\lan x\ran|p|
+|p|^{\rho_1\over\rho_1-\si_1}
+|p|^{\rho_2\over\rho_2-\si_2}\Big\}^{\si_1\over\rho_1}\\
\ns\ds\qq~\le C\Big\{\lan x\ran{}^{\si_1\m\over\rho_1}+\lan
x\ran{}^{\si_1\over\rho_1}|p|^{\si_1\over\rho_1}+|p|^{\si_1\over\rho_1-\si_1}
+|p|^{\si_1\rho_2\over\rho_1(\rho_2-\si_2)}\Big\}.\ea$$
Also, by (\ref{bar u2}), one has
\bel{}\ba{ll}
\ns\ds|\bar u_2|^{\rho_2}\le C\Big\{\lan x\ran{}^\m+\lan
x\ran|p|+|p|\,|u_1|^{\si_1}+|u_1|^{\rho_1}+|p|^{\rho_2\over\rho_2-\si_2}
\Big\}\\
\ns\ds\le C\Big\{\lan x\ran{}^\m+\lan x\ran|p|+
\[\lan x\ran{}^{\si_1\m\over\rho_1}+\lan
x\ran{}^{\si_1\over\rho_1}|p|^{\si_1\over\rho_1}+|p|^{\si_1\over\rho_1-\si_1}+
|p|^{\si_1\rho_2\over\rho_1(\rho_2-\si_2)}\]|p|\\
\ns\ds\qq+\lan x\ran{}^\m+\lan
x\ran|p|+|p|^{\rho_1\over\rho_1-\si_1}+|p|^{\rho_2\over\rho_2-\si_2}\Big\}\\
\ns\ds\le C\Big\{\lan x\ran{}^\m+\lan x\ran|p|+\lan
x\ran{}^{\si_1\m\over\rho_1}|p|+\lan
x\ran{}^{\si_1\over\rho_1}|p|^{\si_1+\rho_1\over\rho_1}+
|p|^{\rho_1\over\rho_1-\si_1}+|p|^{\rho_2\over\rho_2-\si_2}+|p|^{{\si_1\rho_2\over\rho_1(\rho_2-\si_2)}+1}
\Big\}. \ea\ee
Hence,
$$\ba{ll}
\ns\ds|\bar u_2|^{\si_2}\le C\Big\{\lan x\ran{}^\m+\lan
x\ran|p|+\lan x\ran{}^{\si_1\m\over\rho_1}|p|+\lan
x\ran{}^{\si_1\over\rho_1}|p|^{\si_1+\rho_1\over\rho_1}+
|p|^{\rho_1\over\rho_1-\si_1}+|p|^{\rho_2\over\rho_2-\si_2}
+|p|^{{\si_1\rho_2\over\rho_1(\rho_2-\si_2)}+1}
\Big\}^{\si_2\over\rho_2}\\
\ns\ds\qq~\le C\Big\{\lan x\ran{}^{\si_2\m\over\rho_2}+\lan
x\ran{}^{\si_2\over\rho_2}|p|^{\si_2\over\rho_2}+\lan
x\ran{}^{\si_1\si_2\m\over\rho_1\rho_2}|p|^{\si_2\over\rho_2}+\lan
x\ran{}^{\si_1\si_2\over\rho_1\rho_2}|p|^{\si_2(\si_1+\rho_1)
\over\rho_1\rho_2}\\
\ns\ds\qq\qq\qq\qq+|p|^{\si_2\rho_1\over\rho_2(\rho_1-\si_1)}
+|p|^{\si_2\over\rho_2-\si_2}+|p|^{{\si_1\si_2\over\rho_1(\rho_2-\si_2)}+{\si_2\over\rho_2}}
\Big\}.\ea$$
Consequently, for any $(t,x)\in[0,T]\times\dbR^n$, $p,q\in\dbR^n$
and $u_i\in U_i(|x|,|p|\vee|q|)$ ($i=1,2$), we have (without loss of
generality, let $|q|\le|p|$)
\bel{3.24}\ba{ll}
\ns\ds|\dbH(t,x,p,u_1,u_2)-\dbH(t,x,q,u_1,u_2)|\\
\ns\ds\le|p-q|\,|f(t,x,u_1,u_2)|\le L\big(\lan
x\ran+|u_1|^{\si_1}+|u_2|^{\si_2}\big)|p-q|\\
\ns\ds\le C\Big\{\lan x\ran+\lan x\ran{}^{\si_1\m\over\rho_1}+\lan
x\ran{}^{\si_2\m\over\rho_2}+\lan
x\ran{}^{\si_1\over\rho_1}|p|^{\si_1\over\rho_1}+\lan
x\ran{}^{\si_2\over\rho_2}|p|^{\si_2\over\rho_2}+|p|^{\si_1\over\rho_1-\si_1}
+|p|^{\si_2\over\rho_2-\si_2}\\
\ns\ds\qq+|p|^{\si_1\rho_2\over\rho_1(\rho_2-\si_2)}\1n+\1n|p|^{\si_2\rho_1\over\rho_2(\rho_1-\si_1)}
\1n+\1n\lan
x\ran{}^{\si_1\si_2\m\over\rho_1\rho_2}|p|^{\si_2\over\rho_2}\1n+\1n\lan
x\ran{}^{\si_1\si_2\over\rho_1\rho_2}|p|^{\si_2(\si_1+\rho_1)
\over\rho_1\rho_2}\1n+\1n|p|^{{\si_1\si_2\over\rho_1(\rho_2-\si_2)}+{\si_2\over\rho_2}}
\Big\}|p-q|\\
\ns\ds\equiv C\sum_{i=1}^{12}\lan
x\ran{}^{\l_i}\big(|p|\vee|q|\big)^{\n_i}|p-q|.\ea\ee
Due to the fact that the infimum and supremum in (\ref{3.15}) can be
taken on compact sets, we can prove the continuity of $(t,x)\mapsto
H^+(t,x,p)$. \endpf

\ms

A similar result as above can be proved under some much weaker
conditions. In fact, we can relax (H1)--(H2) to the following.

\ms

{\bf(H1)$^*$} Map $f:[0,T]\times\dbR^n\times U_1\times U_2\to\dbR^n$
is continuous and there are constants $\si_1,\si_2\ge0$ and
$\m_0,\m_1,\m_2\in\dbR$ such that
\bel{|f|}\ba{ll}
\ns\ds|f(t,x,u_1,u_2)|\le L\big(\1n\lan x\ran{}^{\m_0}+\lan
x\ran{}^{\m_1}
|u_1|^{\si_1}+\lan x\ran{}^{\m_2}|u_2|^{\si_2}\big),\\
\ns\ds\qq\qq\qq\qq\qq\forall(t,x,u_1,u_2)\in[0,T]\times\dbR^n\times
U_1\times U_2.\ea\ee

\ms

{\bf(H2)$^*$} Map $g:[0,T]\times\dbR^n\times U_1\times U_2\to\dbR$
is continuous and there exist constants $L,c,\rho_1,\rho_2>0$ and
$\bar\m_0,\bar\m_1,\bar\m_2\in\dbR$ such that
\bel{g}\ba{ll}
\ns\ds c\lan x\ran{}^{\bar\m_1}|u_1|^{\rho_1}-L\big(\1n\lan
x\ran{}^{\bar\m_0}+\lan x\ran{}^{\bar\m_2}|u_2|^{\rho_2}\big)\le g(t,x,u_1,u_2)\\
\ns\ds\le L\big(\1n\lan x\ran{}^{\bar\m_0}+\lan
x\ran{}^{\bar\m_1}|u_1|^{\rho_1}\big)-c\lan
x\ran{}^{\bar\m_2}|u_2|^{\rho_2},\qq\forall(t,x,u_1,u_2)\in
[0,T]\times\dbR^n\times U_1\times U_2.\ea\ee

The following result can be proved in the same way as Proposition
3.1.

\ms

\bf Proposition 3.1$^*$. \sl Under {\rm(H1)$^*$--(H2)$^*$} and
{\rm(H3)}, the upper and lower Hamitonians $H^\pm(\cd\,,\cd\,,\cd)$
are well-defined and continuous. Moreover, there are constants
$C>0$, $\n_i\ge0$, and $\l_i\in\dbR$ $(i=1,2,\cds,k)$ such that
\bel{|H|*}\ba{ll}
\ns\ds-L\lan x\ran{}^{\bar\m_0}-L\lan x\ran{}^{\m_0}|p|-C\lan
x\ran{}^{\m_1\rho_1-\bar\m_1\si_1\over\rho_1-\si_1}|p|^{\rho_1\over\rho_1-\si_1}
\le H^\pm(t,x,p)\\
\le L\lan x\ran{}^{\bar\m_0}+L\lan x\ran{}^{\m_0}|p|+C\lan
x\ran{}^{\m_2\rho_2-\bar\m_2\si_2\over\rho_2-\si_2}|p|^{\rho_2\over\rho_2-\si_2},
\qq\qq\forall(t,x,p)\in[0,T]\times\dbR^n\times\dbR^n,\ea\ee
and
\bel{|H-H|*}\ba{ll}
\ns\ds|H^\pm(t,x,p)-H^\pm(t,x,q)|\le C\sum_{i=1}^k\lan
x\ran{}^{\l_i}\big(|p|\vee|q|\big)^{\n_i}|p-q|,\qq\forall(t,x)\in[0,T]\times\dbR^n,~p,q\in\dbR^n.\ea\ee

\ms

\rm

We point out that different from Proposition 3.1, there are more
terms in (\ref{|H-H|*}) than in (\ref{|H-H|}), and the expressions
of $\l_i$ and $\n_i$ are a little more complicated. In fact, instead
of (\ref{3.24}) we can prove the following: (for notational
simplicity, we let $|q|\le|p|$))
$$\ba{ll}
\ns\ds|\dbH(t,x,p,u_1,u_2)-\dbH(t,x,q,u_1,u_2)|\\
\ns\ds\le|p-q|\,|f(t,x,u_1,u_2)|\le L\big(\lan x\ran{}^{\m_0}+\lan
x\ran{}^{\m_1}|u_1|^{\si_1}+\lan
x\ran{}^{\m_2}|u_2|^{\si_2}\big)|p-q|\\
\ns\ds\le C\Big\{\lan x\ran{}^{\m_0}+\lan
x\ran{}^{\m_1-{\si_1\bar\m_1\over\rho_1}}+\lan
x\ran{}^{\m_1+{\si_1(\bar\m_0-\bar\m_1)\over\rho_1}}+\lan
x\ran{}^{\m_1+{\si_1(\m_0-\bar\m_1)\over\rho_1}}|p|^{\si_1\over\rho_1}\\
\ns\ds\qq+\lan
x\ran{}^{\m_1+{\si_1(\m_1-\bar\m_1)\over\rho_1-\si_1}}
|p|^{\si_1\over\rho_1-\si_1}+\lan
x\ran{}^{\m_1+{\si_1[\rho_2(\m_2-\bar\m_1)-\si_2(\bar\m_2-\bar\m_1)]
\over\rho_1(\rho_2-\si_2)}}
|p|^{\si_1\rho_2\over\rho_1(\rho_2-\si_2)}\\
\ns\ds\qq+\lan x\ran{}^{\m_2-{\si_2\bar\m_2\over\rho_2}}+\lan
x\ran{}^{\m_2+{\si_2(\bar\m_0-\bar\m_2)\over\rho_2}}+\lan
x\ran{}^{\m_2+{\si_2(\m_0-\bar\m_2)\over\rho_2}}|p|^{\si_2\over\rho_2}
+\lan
x\ran{}^{\m_2+{\si_2(\m_1-\bar\m_2)\over\rho_2}-{\si_1\si_2\bar\m_1\over\rho_1
\rho_2}}|p|^{\si_2\over\rho_2}\\
\ns\ds\qq+\lan
x\ran{}^{\m_2+{\si_2(\m_1-\bar\m_2)\over\rho_2}+{\si_1\si_2(\bar\m_0-\bar\m_1)
\over\rho_1\rho_2}}|p|^{\si_2\over\rho_2}+\lan
x\ran{}^{\m_2+{\si_2(\m_1-\bar\m_2)\over\rho_2}+{\si_1\si_2(\m_0-\bar\m_1)\over
\rho_1\rho_2}}|p|^{\si_2(\si_1+\rho_1)\over\rho_1\rho_2}\\
\ns\ds\qq+\lan
x\ran{}^{\m_2+{\si_2(\m_1-\bar\m_2)\over\rho_2}+{\si_1\si_2(\m_1-\bar\m_1)
\over\rho_2(\rho_1-\si_1)}}|p|^{\si_2\rho_1\over\rho_2(\rho_1-\si_1)}\\
\ns\ds\qq+\lan
x\ran{}^{\m_2+{\si_2(\m_1-\bar\m_2)\over\rho_2}+{\si_1\si_2
[\rho_2(\m_2-\bar\m_1)-\si_2(\bar\m_2-\bar\m_1)]
\over\rho_1\rho_2(\rho_2-\si_2)}}
|p|^{{\si_1\si_2\over\rho_1(\rho_2-\si_2)}+{\si_2\over\rho_2}}\\
\ns\ds\qq+\lan
x\ran{}^{\m_2+{\si_2(\m_1\rho_1-\bar\m_1\si_1)\over\rho_2(\rho_1-\si_1)}
-{\si_2\bar\m_2\over\rho_2}}|p|^{\si_2\rho_1\over\rho_2(\rho_1-\si_1)}
+\lan
x\ran{}^{\m_2+{\si_2(\m_2\rho_2-\bar\m_2\si_2)\over\rho_2(\rho_2-\si_2)}
-{\si_2\bar\m_2\over\rho_2}}
|p|^{\si_2\over\rho_2-\si_2}\Big\}|p-q|\\
\ns\ds\equiv C\sum_{i=1}^{16}\lan
x\ran{}^{\l_i}\big(|p|\vee|q|\big)^{\n_i}|p-q|.\ea$$
Note that (\ref{3.24}) is a special case of the above with:
$$\m_0=1,~\bar\m_0=\m,~\m_1=\m_2=\bar\m_1=\bar\m_2=0.$$

\ms

\section{Uniqueness of Viscosity Solution}

\ms

Consider the following HJ inequalities:
\bel{HJ-sub}\left\{\ba{ll}
\ns\ds V_t(t,x)+H(t,x,V_x(t,x))\ge0,\qq(t,x)\in[0,T]\times\dbR^n,\\
\ns\ds V(T,x)\le h(x),\qq\qq x\in\dbR^n,\ea\right.\ee
and
\bel{HJ-super}\left\{\ba{ll}
\ns\ds V_t(t,x)+H(t,x,V_x(t,x))\le0,\qq(t,x)\in[0,T]\times\dbR^n,\\
\ns\ds V(T,x)\ge h(x),\qq\qq x\in\dbR^n,\ea\right.\ee
as well as the following HJ equation:
\bel{HJ}\left\{\ba{ll}
\ns\ds V_t(t,x)+H(t,x,V_x(t,x))=0,\qq(t,x)\in[0,T]\times\dbR^n,\\
\ns\ds V(T,x)=h(x),\qq\qq x\in\dbR^n.\ea\right.\ee
We recall the following definition.

\ms

\bf Definition 4.1. \rm (i) A continuous function $V(\cd\,,\cd)$ is
called a {\it viscosity sub-solution} of (\ref{HJ-sub}) if
$$V(T,x)\le h(x),\qq\forall x\in\dbR^n,$$
and for any continuous differentiable function $\f(\cd\,,\cd)$, if
$(t_0,x_0)\in[0,T)\times\dbR^n$ is a local maximum of $(t,x)\mapsto
V(t,x)-\f(t,x)$, then
$$\f_t(t_0,x_0)+H(t_0,x_0,\f_x(t_0,x_0))\ge0.$$

(ii) A continuous function $V(\cd\,,\cd)$ is called a {\it viscosity
super-solution} of (\ref{HJ-super}) if
$$V(T,x)\ge h(x),\qq\forall x\in\dbR^n,$$
and for any continuous differentiable function $\f(\cd\,,\cd)$, if
$(t_0,x_0)\in[0,T)\times\dbR^n$ is a local minimum of $(t,x)\mapsto
V(t,x)-\f(t,x)$, then
$$\f_t(t_0,x_0)+H(t_0,x_0,\f_x(t_0,x_0))\le0.$$

(iii) A continuous function $V(\cd\,,\cd)$ is called a {\it
viscosity solution} of (\ref{HJ}) if it is a viscosity sub-solution
of (\ref{HJ-sub}) and a viscosity super-solution of
(\ref{HJ-super}).

\ms

The following lemma is taken from \cite{Crandall-Lions 1986}.

\ms

\bf Lemma 4.2. \sl Suppose $H:[0,T]\times\dbR^n\times\dbR^n\to\dbR$
is continuous. Let $V(\cd\,,\cd)$ and $\h V(\cd\,,\cd)$ be a
viscosity sub- and super-solutions of $(\ref{HJ-sub})$ and
$(\ref{HJ-super})$, respectively. Then
$$W(t,x,y)=V(t,x)-\h V(t,y),\qq(t,x,y)\in[0,T]\times\dbR^n\times\dbR^n$$
is a viscosity sub-solution of the following:
$$\left\{\ba{ll}
\ns\ds
W_t(t,x,y)+H(t,x,W_x(t,x,y))-H(t,y,-W_x(t,x,y))\ge0,\qq(t,x,y)\in[0,T]\times\dbR^n\times\dbR^n,\\
\ns\ds W(T,x,y)\le0,\qq\qq(x,y)\in\dbR^n\times\dbR^n.\ea\right.$$

\ms

\rm

Now for HJ equation (\ref{HJ}), we assume the following.

\ms

{\bf(HJ)} The maps $H:[0,T]\times\dbR^n\times\dbR^n\to\dbR$ and
$h:\dbR^n\to\dbR$ are continuous and there are constants $K_0>0$,
$\m\ge1$, and $\l_i,\n_i\ge0$ $(i=1,2,\cds,k$) with
\bel{ml}\l_i+(\m-1)\n_i\le1,\qq1\le i\le k,\ee
and a continuous function
$\o:[0,\infty)^3\to[0,\infty)$ with property $\o(r,s,0)=0$, such
that
\bel{H<}\ba{ll}
\ns\ds|H(t,x,p)-H(t,y,p)|\le\o\big(|x|+|y|,|p|,|x-y|\big),\qq\forall
t\in[0,T],~x,y,p\in\dbR^n,\\
\ns\ds|H(t,x,p)-H(t,x,q)|\le K_0\sum_{i=1}^k\lan
x\ran{}^{\l_i}\big(|p|\vee|q|\big)^{\n_i}|p-q|,\qq\forall
t\in[0,T],~x,p,q\in\dbR^n.\ea\ee
and
\bel{h<}|h(x)-h(y)|\le K_0\big(\1n\lan x\ran\vee\lan
y\ran\big)^{\m-1}|x-y|,\qq\forall x,y\in\dbR^n.\ee

Our main result of this section is the following.

\ms

\bf Theorem 4.3. \sl Let {\rm(HJ)} hold. Suppose $V(\cd\,,\cd)$ and
$\h V(\cd\,,\cd)$ are the viscosity sub- and super-solution of
$(\ref{HJ-sub})$ and $(\ref{HJ-super})$, respectively. Moreover, let
\bel{|V-V|}|V(t,x)-V(t,y)|,|\h V(t,x)-\h V(t,y)|\le K\big(\1n\lan
x\ran\vee\lan y\ran\big)^{\m-1}|x-y|,\q\forall
t\in[0,T],~x,y\in\dbR^n,\ee
for some $K>0$. Then
\bel{V<V}V(t,x)\le\h V(t,x),\qq\forall(t,x)\in[0,T]\times\dbR^n.\ee

\ms

\rm

A similar result as above was proved in \cite{Crandall-Lions 1987},
with most technical details omitted. Our conditions are a little
different from those assumed in \cite{Crandall-Lions 1987}. For
readers' convenience, we provide a detailed proof here.

\ms

\it Proof. \rm Suppose $(\bar t,\bar x)\in[0,T)\times\dbR^n$ such
that
$$V(\bar t,\bar x)-\h V(\bar t,\bar x)>0.$$
Let $C_0,\b>0$ be undetermined. Define
$$Q\equiv Q(C_0,\b)=\Big\{(t,x)\in[0,T]\times\dbR^n\bigm|\lan
x\ran\le\lan\bar x\ran e^{C_0(t-\bar t)+\b}\Big\},$$
and
$$G\equiv G(C_0,\b)=\Big\{(t,x,y)\in[0,T]\times\dbR^n\times\dbR^n\bigm|(t,x),(t,y)\in Q\Big\}.$$
Now, for $\d>0$ small, define
$$\ba{ll}
\ns\ds\psi(t,x)\equiv\psi^{C_0,\d}(t,x)=\[{\lan x\ran\over\lan\bar
x\ran}e^{C_0(\bar t-t)}\]^{1\over\d}\equiv
e^{{1\over\d}\big[\log{\lan x\ran\over\lan\bar x\ran}+C_0(\bar
t-t)\big]}.\ea$$
Then
$$\psi(\bar t,\bar x)=1,\qq\psi(T,x)=\[{\lan x\ran\over\lan\bar x\ran}
e^{C_0(\bar t-T)}\]^{1\over\d},$$
and
$$\psi_t(t,x)=-{C_0\psi(t,x)\over\d}\,,\qq\psi_x(t,x)={x\psi(t,x)\over\d\lan x\ran{}^2}\,.$$
For any $(t,x)\in\bar Q$, we have
$$\lan x\ran\le\lan\bar x\ran e^{\b+C_0(t-\bar t)}
\le\lan\bar x\ran e^{\b+C_0(T-\bar t)}.$$
Thus, $Q$ is bounded and $\bar G$ is compact. We introduce
$$\Psi(t,x,y)=V(t,x)-\h
V(t,y)-{|x-y|^2\over\e}-\si\psi(t,x)-{\si(T-t)\over T-\bar
t}\,,\qq(t,x,y)\in\bar G,$$
where $\e>0$ small and
$$0<\si\le{V(\bar t,\bar x)-\h V(\bar t,\bar x)\over3}.$$
Clearly,
\bel{psi>si}\Psi(\bar t,\bar x,\bar x)=V(\bar t,\bar x)-\h V(\bar
t,\bar x)-\si-\si\ge3\si-2\si=\si>0.\ee
Since $\Psi(\cd\,,\cd\,,\cd)$ is continuous on the compact set $\bar
G$, we may let $(t_0,x_0,y_0)\in\bar G$ be a maximum of
$\Psi(\cd\,,\cd\,,\cd)$ over $\bar G$. By the optimality of
$(t_0,x_0,y_0)$, we have
$$\ba{ll}
\ns\ds V(t_0,x_0)-\h
V(t_0,x_0)-\si\psi(t_0,x_0)-{\si(T-t_0)\over T-\bar t}=\Psi(t_0,x_0,x_0)\\
\ns\ds\le\Psi(t_0,x_0,y_0)=V(t_0,x_0)-\h
V(t_0,y_0)-{|x_0-y_0|^2\over\e}-\si\psi(t_0,x_0)-{\si(T-t_0)\over
T-\bar t},\ea$$
which implies
$${|x_0-y_0|^2\over\e}\le\h V(t_0,x_0)-\h V(t_0,y_0)\le
K\big(\1n\lan x_0\ran\vee\lan y_0\ran\big)^{\m-1}|x_0-y_0|.$$
Thus,
$$|x_0-y_0|\le K\big(\1n\lan x_0\ran\vee\lan y_0\ran\big)^{\m-1}\e.$$
Now, if $t_0=T$, then
$$\lan x_0\ran,~\lan y_0\ran\le\lan\bar x\ran e^{\b+C_0(T-\bar t)}.$$
Hence,
$$\ba{ll}
\ns\ds\Psi(T,x_0,y_0)=h(x_0)-h(y_0)-{|x_0-y_0|^2\over\e}-\si\psi(T,x_0)\le
K_0\big(\1n\lan x_0\ran\vee\lan y_0\ran\big)^{\m-1}|x_0-y_0|\\
\ns\ds\qq\qq\q\;\le K_0K\big(\1n\lan x_0\ran\vee\lan
y_0\ran\big)^{2(\m-1)}\e\le K_0K\lan\bar
x\ran{}^{2(\m-1)}e^{2(\m-1)[\b+C_0(T-\bar t)]}\e.\ea$$
Thus, for $\e>0$ small enough, the following holds:
$$\Psi(T,x_0,y_0)<\si\le\Psi(\bar t,\bar x,\bar x)\le\Psi(t_0,x_0,y_0),$$
which means that $t_0\in[0,T)$. Next, we note that for
$(t,x)\in\big(\pa Q\big)\cap\big[(0,T)\times\dbR^n\big]$, one has
$$\log{\lan x\ran\over\lan\bar x\ran}+C_0(\bar t-t)=\b,\q\hb{ and }\q0<t<T,$$
which implies
\bel{psi-infty}\psi(t,x)=e^{\b\over\d}\to\infty,\qq\d\to0,\q\hb{uniformly
in }(t,x)\in\big(\pa Q\big)\cap\big[(0,T)\times\dbR^n\big].\ee
This implies that for $\d>0$ small (only depending on $\b$),
$$(t_0,x_0,y_0)\in G\cup\[\{0\}\times\dbR^n\times\dbR^n\].$$
By Lemma 4.2, we have
$$\ba{ll}
\ns\ds0\le\si\psi_t(t_0,x_0)-{\si\over T-\bar
t}+H\(t_0,x_0,{2(x_0-y_0)\over\e}+\si\psi_x(t_0,x_0)\)
-H\(t_0,y_0,-{2(y_0-x_0)\over\e}\)\\
\ns\ds=\si\psi_t(t_0,x_0)-{\si\over T-\bar
t}+H\(t_0,x_0,{2(x_0-y_0)\over\e}+\si\psi_x(t_0,x_0)\)
-H\(t_0,x_0,{2(x_0-y_0)\over\e}\)\\
\ns\ds\qq\qq+H\(t_0,x_0,{2(x_0-y_0)\over\e}\)-H\(t_0,y_0,{2(x_0-y_0)\over\e}\)\\
\ns\ds\le\1n\si\psi_t(t_0,x_0)\1n-\1n{\si\over T\1n-\bar t}+\1n
K_0\sum_{i=1}^k\lan x_0\ran{}^{\l_i}\({2|x_0\1n-\1n y_0|\over\e}\1n
+\1n\si|\psi_x(t_0,x_0)|\)^{\n_i}\si|\psi_x(t_0,x_0)|\\
\ns\ds\qq+\o\(|x_0|+|y_0|,{2|x_0-y_0|\over\e},|x_0-y_0|\)\\
\ns\ds\le-\si{C_0\over\d}\psi(t_0,x_0)-{\si\over T-\bar t}+\si
K_0\sum_{i=1}^k\lan x_0\ran{}^{\l_i}\(2K\big(\1n\lan x_0\ran\vee\lan
y_0\ran\big)^{\m-1}+{\si\psi(t_0,x_0)\over\d\lan
x_0\ran}\)^{\n_i}\]{\psi(t_0,x_0)\over\d\lan x_0\ran}\\
\ns\ds\qq+\o\(|x_0|+|y_0|,{2|x_0-y_0|\over\e},|x_0-y_0|\).\ea$$
Note that $(t_0,x_0,y_0)\equiv(t_{0,\e},x_{0,\e},y_{0,\e})\in\bar
G(C_0,\b)$ (a fixed compact set). Let $\e\to0$ along a suitable
sequence, we have $|x_{0,\e}-y_{0,\e}|\to0$. For notational
simplicity, we denote
$(t_{0,\e},x_{0,\e},y_{0,\e})\to(t_0,x_0,x_0)$. In the above, by
canceling $\si$, and then send $\e\to0$ and $\si\to0$, one obtains
(canceling $\si$)
$$\ba{ll}
\ns\ds{1\over T-\bar
t}\le-{C_0\psi(t_0,x_0)\over\d}+K_0\sum_{i=1}^k\lan
x_0\ran{}^{\l_i}\(2K\1n\lan
x_0\ran{}^{\m-1}\)^{\n_i}{\psi(t_0,x_0)\over\d\lan x_0\ran}\\
\ns\ds=-\Big\{C_0-K_0\sum_{i=1}^k(2K)^{\n_i}\lan
x_0\ran{}^{\l_i+(\m-1)\n_i-1}\Big\}{\psi(t_0,x_0)\over\d}\\
\ns\ds\le-\Big\{C_0-K_0\sum_{i=1}^k(2K)^{\n_i}\Big\}{\psi(t_0,x_0)\over\d}\equiv-\big(C_0-\wt
K_0\big){\psi(t_0,x_0)\over\d}\,.\ea$$
Thus, by taking $C_0>\wt K_0$, we obtain a contradiction, proving
our conclusion. \endpf

\ms

\rm

We now make some comments on the uniqueness/non-uniqueness of
viscosity solutions. First of all, let us look at the following
example which is adopted from \cite{Crandall-Lions 1983, Biton
2001},

\ms

\bf Example 4.4. \rm It is known that there are two different
bounded strictly increasing continuous differentiable functions
$f_i:\dbR\to\dbR$ ($i=1,2$) such that
$$b(x)\equiv f_1'(f_1^{-1}(x))=f_2'(f_2^{-1}(x)),\qq x\in\dbR.$$
Further, if we define
$$X^i(t;x_0)=f_i(t+f_i^{-1}(x_0)),\qq t\in\dbR$$
then $X^1(\cd\,;x_0)$ and $X^2(\cd\,;x_0)$ are two different
solutions to the following initial value problem:
$$\left\{\ba{ll}
\ns\ds{d\over dt}X(t;x_0)=b\big(X(t;x_0)\big),\qq t\in\dbR,\\
\ns\ds X(0;x_0)=x_0.\ea\right.$$
By defining
$$V^i(t,x)=h(X^i(T-t;x)),\qq(t,x)\in[0,T]\times\dbR,$$
we obtain two different viscosity solutions to the following HJ
equation:
\bel{HJb}\left\{\ba{ll}
\ns\ds V_t(t,x)+b(x)V_x(t,x)=0,\qq(t,x)\in[0,T]\times\dbR,\\
\ns\ds V(T,x)=h(x),\qq x\in\dbR.\ea\right.\ee
Therefore, the viscosity solution to the above HJ equation is not
unique in the set of continuous functions. However, we note that in
the current case,
$$H(t,x,p)=b(x)p,\qq(t,x,p)\in[0,T]\times\dbR\times\dbR.$$
Thus,
$$|H(t,x,p)-H(t,x,q)|\le C|p-q|,\qq\forall t\in[0,T],~x,p,q\in\dbR,$$
which means that (\ref{H<}) holds with $k=1$, $\l_1=\n_1=0$. Hence,
for any $\m\ge1$, as long as (\ref{h<}) holds, viscosity solution to
(\ref{HJb}) is unique in the class of continuous functions
satisfying (\ref{|V-V|}).

\ms

\bf Example 4.5. \rm Consider
$$-x^2-axV_x(x)+|V_x(x)|^2=0,\qq x\in\dbR,$$
with $a\ge0$. Thus,
$$H(x,p)=-x^2-axp+p^2,\qq(x,p)\in\dbR^2.$$
Then let
$$V(x)=\l x^2,\qq x\in\dbR.$$
We should have
$$0=-1-2a\l+4\l^2.$$
Hence,
$$\l={2a\pm\sqrt{4a^2+16}\over8}={a\pm\sqrt{a^2+4}\over4}.$$
Therefore, there are two solutions to the HJ equation:
$$V^\pm(x)={a\pm\sqrt{a^2+4\,}\over4}\,x^2,\qq x\in\dbR.$$
Both of these solutions are analytic. Note that
$$|H(x,p)-H(x,q)|\le\(a|x|+2\big(|p|\vee|q|\big)\)|p-q|,\qq
x,p,q\in\dbR,$$
and
$$|V^\pm(x)-V^\pm(y)|\le{|a\pm\sqrt{a^2+4\,}|\over4}
\big(\1n\lan x\ran\vee\lan y\ran\1n\big)|x-y|.$$
Thus, in our terminology, $\m=2$, $k=2$ with
$$\l_1=0,\q\n_1=1,\q\l_2=0,\q\n_2=1.$$
Consequently,
$$\l_i+(\m-1)\n_i=1,\qq i=1,2.$$
This means that although (\ref{ml}) is satisfied, the corresponding
HJ equation has more than one viscosity solution. This example shows
that stationary problems are different from evolution problems, as
far as the uniqueness of viscosity solution is concerned.

\section{Upper and Lower Value Functions}

In this section, we are going to define the upper and lower value
functions via the so-called Elliott--Kalton strategies. Some basic
properties of upper and lower value functions will be established
carefully.

\subsection{State trajectories and Elliott--Kalton strategies}

Let us introduce the following hypotheses which are strengthened
versions of (H1)--(H3).

\ms

{\bf(H1)$'$} Map $f:[0,T]\times\dbR^n\times U_1\times U_2\to\dbR^n$
satisfies (H1). Moreover, for some $\m_0,\m_1,\m_2$,
\bel{|f-f|}\ba{ll}
\ns\ds|f(t,x,u_1,u_2)-f(t,y,u_1,u_2)|\\
\ns\ds\le\big[\big(\1n\lan x\ran\vee\lan
y\ran\big)^{\m_0}+\big(\1n\lan x\ran\vee\lan
y\ran\big)^{\m_1}|u_1|^{\si_1}+\big(\1n\lan x\ran\vee\lan
y\ran\big)^{\m_2}|u_2|^{\si_2}\big]|x-y|,\\
\ns\ds\qq\qq\qq\qq\qq\qq\forall(t,u_1,u_2)\in[0,T]\times U_1\times
U_2,~x,y\in\dbR^n,\ea\ee
and
\bel{f-f}\ba{ll}
\lan f(t,x,u_1,u_2)-f(t,y,u_1,u_2),x-y\ran\le L|x-y|^2,\\
\ns\ds\qq\qq\qq\qq\qq\qq\qq\forall(t,u_1,u_2)\in[0,T]\times
U_1\times U_2,~x,y\in\dbR^n.\ea\ee

\ms

We note that condition (\ref{|f-f|}) implies the local Lipschitz
continuity of the map $x\mapsto f(t,x,u_1,u_2)$, with the Lipschtiz
constant possibly depending on $|u_1|^{\si_1}$ and $|u_2|^{\si_2}$.
This is the case if we are considering AQ two-person zero-sum
differential games (see Section 2). On the other hand, condition
(\ref{f-f}) will be used to establish the local Lipschitz continuity
of the upper and lower value functions, with the Lipschitz constant
being of polynomial order of $\lan x\ran\vee\lan y\ran$. It is
important that the right hand side of (\ref{f-f}) is independent of
$(u_1,u_2)$; Otherwise, the Lipschitz constant of the upper and
lower value functions will be some exponential function of $\lan
x\ran\vee\lan y\ran$, for which we do not know if the uniqueness of
viscosity solution to the corresponding HJI equation holds. By the
way, we point out that (\ref{f-f}) does not imply the local
Lipschitz continuity of the map $x\mapsto f(t,x,u_1,u_2)$. For
example, $f(x)=x^{1\over3}$, with $x\in\dbR$.

\ms

{\bf(H2)$'$} Map $g:[0,T]\times\dbR^n\times U_1\times U_2\to\dbR$
satisfies (H2). Moreover,
\bel{g-g}\ba{ll}
\ns\ds|g(t,x,u_1,u_2)-g(t,y,u_1,u_2)|\le\big[\big(\1n\lan
x\ran\vee\lan
y\ran\big)^{\m-1}+|u_1|^{\rho_1(\m-1)\over\m}+|u_2|^{\rho_2(\m-1)\over\m}\big]|x-y|,\\
\ns\ds\qq\qq\qq\qq\qq\qq\qq\forall(t,u_1,u_2)\in[0,T]\times
U_1\times U_2,~x,y\in\dbR^n.\ea\ee
Also, map $h:\dbR^n\to\dbR$ is continuous and
\bel{h}\left\{\ba{ll}
\ns\ds|h(x)-h(y)|\le L\big(\1n\lan x\ran\vee\lan
y\ran\big)^{\m-1}|x-y|,\qq\forall x,y\in\dbR^n,\\
\ns\ds|h(0)|\le L.\ea \right.\ee

\ms

Further, the compatibility hypothesis (H3) is now replaced by the
following:

\ms

{\bf(H3)$'$} The constants $\si_1,\si_2,\rho_1,\rho_2,\m$ appear in
(H1)$'$--(H2)$'$ satisfy the following:
\bel{5.4}\si_i\m<\rho_i,\qq i=1,2.\ee

\ms

Let us first present the following Gronwall type inequality.

\ms

\bf Lemma 5.1. \sl Let $\th,\a,\b:[t,T]\to\dbR_+$ and $\th_0\ge0$
satisfy
\bel{}\th(s)^2\le\th_0^2+\int_t^s\big[\a(r)\th(r)^2+\b(r)\th(r)\big]dr,\qq
s\in[t,T].\ee
Then
\bel{Gronwall}\th(s)\le
e^{{1\over2}\int_t^T\a(\t)d\t}\th_0+{1\over2}e^{\int_t^T\a(\t)d\t}\int_t^s\b(r)dr,
\qq s\in[t,T].\ee

\ms

\it Proof. \rm First, by the usual Gronwall's inequality, we have
$$\ba{ll}
\ns\ds\th(s)^2\le
e^{\int_t^s\a(\t)d\t}\th_0^2+\int_t^se^{\int_r^s\a(\t)d\t}\b(r)\th(r)dr\\
\ns\ds\qq~\le
e^{\int_t^T\a(\t)d\t}\th_0^2+e^{\int_t^T\a(\t)d\t}\int_t^s\b(r)\th(r)dr\equiv\Th(s).\ea$$
Then
$${d\over
ds}\sqrt{\Th(s)}={1\over2}\Th(s)^{-{1\over2}}\dot\Th(s)={1\over2}\Th(s)^{-{1\over2}}e^{\int_t^T\a(\t)d\t}
\b(s)\th(s)\le{1\over2}e^{\int_t^T\a(\t)d\t}\b(s).$$
Consequently,
$$\th(s)\le\sqrt{\Th(s)}\le e^{{1\over2}\int_t^T\a(\t)d\t}\th_0+{1\over2}e^{\int_t^T\a(\t)d\t}\int_t^s\b(r)dr,
\qq s\in[t,T],$$
proving our conclusion. \endpf

\ms

We now prove the following result concerning the state trajectories.

\ms

\bf Proposition 5.2. \sl Let {\rm(H1)$'$} hold. Then, for any
$(t,x)\in[0,T)\times\dbR^n$,
$(u_1(\cd),u_2(\cd))\in\cU_1^{\si_1}[t,T]\times\cU_2^{\si_2}[t,T]$,
state equation $(\ref{1.1})$ admits a unique solution
$y(\cd)\1n\equiv\1n y(\cd\,;t,x,u_1(\cd),u_2(\cd))\1n\equiv
y_{t,x}(\cd)$. Moreover, there exists a constant $C_0>0$ only
depends on $L,T,t$ such that
\bel{|y(s)|}\lan y_{t,x}(s)\ran\le C_0\Big\{\lan
x\ran+\int_t^s\big(|u_1(r)|^{\si_1} +|u_2(r)|^{\si_2}\big)dr\Big\},
\qq\qq s\in[t,T],\ee
\bel{y-x}|y_{t,x}(s)-x|\le C_0\Big\{\2n\lan
x\ran(s-t)+\int_t^s\big(|u_1(r)|^{\si_1}+|u_2(r)|^{\si_2}\big)dr\Big\},\q
s\in[t,T].\ee
and for $(\bar t,\bar x)\in[0,T]\times\dbR^n$ with $\bar t\in[t,T]$,
and $y_{\bar t,\bar x}(\cd)\equiv y(\cd\,;\bar t,\bar
x,u_1(\cd),u_2(\cd))$
\bel{}\ba{ll}
\ns\ds|y_{t,x}(s)-y_{\bar t,\bar x}(s)|\le C_0\Big\{|x-\bar x|+\lan
x\ran(\bar t-t)+\int_t^{\bar
t}\big(|u_1(r)|^{\si_1}+|u_2(r)|^{\si_2}\big)dr\Big\},\qq
s\in[t,T].\ea\ee

\ms

\it Proof. \rm First, under (H1)$'$, for any
$(t,x)\in[0,T)\times\dbR^n$, and any
$(u_1(\cd),u_2(\cd))\in\cU_1^{\si_1}[t,T]\times\cU_2^{\si_2}[t,T]$,
the map $y\mapsto f(s,y,u_1(s),u_2(s))$ is locally Lipschitz
continuous. Thus, state equation (\ref{1.1}) admits a unique local
solution $y(\cd)=y(\cd\,;t,x,u_1(\cd),u_2(\cd))$. Next, by
(\ref{f-f}), we have
$$\ba{ll}
\ns\ds\lan x,f(t,x,u_1,u_2)\ran=\lan
x,f(t,x,u_1,u_2)-f(t,0,u_1,u_2)\ran+\lan x,f(t,0,u_1,u_2)\ran\\
\ns\ds\le
L|x|^2+L|x|\big(\1n1+|u_1|^{\si_1}+|u_2|^{\si_2}\big),\qq\qq\forall(t,x,u_1,u_2)\in[0,T]\times\dbR^n\times
U_1\times U_2.\ea$$
Thus,
$$\ba{ll}
\ns\ds\lan y(s)\ran{}^2=\lan x\ran{}^2+2\int_t^s\lan y(r),f(r,y(r),u_1(r),u_2(r))\ran dr\\
\ns\ds\qq\qq\le\lan x\ran{}^2+2\int_t^sL\(\lan y(r)\ran{}^2+\lan
y(r)\ran\big(1+|u_1(r)|^{\si_1}+|u_2(r)|^{\si_2}\big)\)dr.\ea$$
%
%
Then, it follows from Lemma 5.1 that
$$\lan y(s)\ran\le e^{L(T-t)}\lan
x\ran+Le^{2L(T-t)}\int_t^s\big(1+|u_1(r)|^{\si_1}
+|u_2(r)|^{\si_2}\big)dr.$$
This implies that the solution $y(\cd)$ of the state equation
(\ref{1.1}) globally exists on $[t,T]$ and (\ref{|y(s)|}) holds.
Also, we have
$$\ba{ll}
\ns\ds|y(s)-x|^2=2\int_t^s\lan y(r)-x,f(r,y(r),u_1(r),u_2(r))\ran dr\\
\ns\ds\qq\qq\q\le2\int_t^s\(L|y(r)-x|^2+\lan
y(r)-x,f(r,x,u_1(r),u_2(r))\ran\)dr\\
\ns\ds\qq\qq\q\le2L\int_t^s\(|y(r)-x|^2+|y(r)-x|\big(\1n\lan
x\ran+|u_1(r)|^{\si_1}+|u_2(r)|^{\si_2}\big)\)dr.\ea$$
Thus, by Lemma 5.2 again, we obtain (\ref{y-x}).

\ms

Now, for any $(t,x),(\bar t,\bar x)\in[0,T]\times\dbR^n$, with $0\le
t\le\bar t<T$, denote $y_{t,x}(\cd)=y(\cd\,;t,x,u_1(\cd),u_2(\cd))$,
and $y_{\bar t,\bar x}(\cd)=y(\cd\,;\bar t,\bar
x,u_1(\cd),u_2(\cd))$. Then for $s\in[\bar t,T]$, we have
$$\ba{ll}
\ns\ds|y_{t,x}(s)-y_{\bar t,\bar x}(s)|^2=|y_{t,x}(\bar t)-\bar
x|^2\\
\ns\ds\qq\qq\qq\qq\q+2\int_{\bar t}^s\lan y_{t,x}(r)-y_{\bar t,\bar
x}(r),f(r,y_{t,x}(r),u_1(r),u_2(r))-f(r,y_{\bar t,\bar
x}(r),u_1(r),u_2(r))\ran dr\\
\ns\ds\qq\qq\qq\qq\le|y_{t,x}(\bar t)-x|^2+2L\int_{\bar
t}^s|y_{t,x}(r)-y_{\bar t,\bar x}(r)|^2dr.\ea$$
Thus, it follows from the Gronwall's inequality that
$$\ba{ll}
\ns\ds|y_{t,x}(s)-y_{\bar t,\bar x}(s)|\le e^{L(s-\bar
t)}|y_{t,x}(\bar t)-x|\le e^{L(s-\bar t)}\big(|x-\bar
x|+|y_{t,x}(\bar t)-x|\big)\\
\ns\ds\le e^{L(s-\bar t)}\Big\{|x-\bar x|+Le^{2L(T-t)}\(\2n\lan
x\ran(\bar t-t)+\int_t^{\bar t}\2n|u_1(r)|^{\si_1}dr+\int_t^{\bar
t}\2n|u_2(r)|^{\si_2}dr\)\Big\}\\
\ns\ds\le C\Big\{|x-\bar x|+\lan x\ran(\bar t-t)+\int_t^{\bar
t}\big(|u_1(r)|^{\si_1}+|u_2(r)|^{\si_2}\big)dr\Big\}.\ea$$
This completes the proof. \endpf

\ms

From the above proposition, together with (H2)$'$, we see that for
any $u_i(\cd)\in\cU_i^{\rho_i}[t,T]$ (which is smaller than
$\cU_i^{\si_i}[t,T]$), $i=1,2$, the performance functional
$J(t,x;u_1(\cd),u_2(\cd))$ is well-defined. Let us now introduce the
following definition which is a modification of the notion
introduced in \cite{Elliott-Kalton 1972}.

\ms

\bf Definition 5.3. \rm A map
$\a_1:\cU_2^1[t,T]\to\cU_1^\infty[t,T]$ is called an Elliott--Kalton
(E-K, for short) strategy for Player 1 if it is {\it
non-anticipating}, namely, for any $u_2(\cd),\bar
u_2(\cd)\in\cU_2^1[t,T]$, and any $\hat t\in[t,T]$,
$$\a_1[u_2(\cd)](s)=\a_1[\bar u_2(\cd)](s),\qq\ae s\in[t,\hat t],$$
provided
$$u_1(s)=\bar u_1(s),\qq\ae s\in[t,\hat t].$$
The set of all E-K strategies for Player 1 is denoted by
$\cA_1[t,T]$. An E-K strategy
$\a_2:\cU_2^1[t,T]\to\cU_1^\infty[t,T]$ for Player 2 can be defined
similarly. The set of all E-K strategies for Player 2 is denoted by
$\cA_2[t,T]$.

\ms

Note that as far as the state equation is concerned, one could
define an E-K strategy $\a_1$ for Player I as a map
$\a_1:\cU^{\si_2}_2[t,T]\to\cU^{\si_1}_1[t,T]$. Whereas, as far as
the performance functional is concerned, one might have to
restrictively define
$\a_1:\cU_2^{\rho_2}[t,T]\to\cU_1^{\rho_1}[t,T]$. We note that the
numbers $\si_1,\si_2,\rho_1,\rho_2$ appeared in (H1)$'$--(H2)$'$
might not be the ``optimal'' ones, in some sense (for example,
$\si_1$ and $\si_2$ might be larger than necessary, and $\rho_1$ and
$\rho_2$ could be smaller than they should be, and so on). Our above
definition is somehow ``universal''. The domain $\cU_2^1[t,T]$ of
$\a_1$ is large enough to cover possible $u_2(\cd)$ in some larger
space than $\cU_2^{\si_2}[t,T]$, and the co-domain
$\cU_1^\infty[t,T]$ is large enough so that the integrability of
$\a_1[u_2(\cd)]$ is ensured and the supremum  will remain the same
due to the density of $\cU_1^\infty[t,T]$ in $\cU_1^{\rho_1}[t,T]$.
In what follows, we simply denote
$$\cU_i[t,T]=\cU_i^\infty[t,T],\qq i=1,2.$$
Recall that $0\in U_i$ ($i=1,2$). For later convenience, we
hereafter let $u^0_1(\cd)\in\cU_1[t,T]$ and
$u^0_2(\cd)\in\cU_2[t,T]$ be defined by
$$u^0_1(s)=0,\q u^0_2(s)=0,\qq\forall s\in[t,T],$$
and let $\a_1^0\in\cA_1[t,T]$ be the E-K strategy that
$$\a_1^0[u_2(\cd)](s)=0,\qq\forall s\in[t,T],\q
u_2(\cd)\in\cU_2^1[t,T].$$
We call such an $\a_1^0$ the {\it zero E-K strategy} for Player 1.
Similarly, we define zero E-K strategy $\a_2^0\in\cA_2[t,T]$ for
Player 2.

\ms

Now, we define
\bel{upper-lower}\left\{\ba{ll}
\ns\ds
V^+(t,x)=\sup_{\a_2\in\cA_2[t,T]}\inf_{u_1(\cd)\in\cU_1[t,T]}J(t,x;u_1(\cd),\a_2[u_1(\cd)]),\\
\ns\ds
V^-(t,x)=\inf_{\a_1\in\cA_1[t,T]}\sup_{u_2(\cd)\in\cU_2[t,T]}J(t,x;\a_1[u_2(\cd)],u_2(\cd)).\ea
\right.\qq(t,x)\in[0,T]\times\dbR^n,\ee
which are called {\it upper} and {\it lower value functions} of our
two-person zero-sum differential game.

\ms

\subsection{Upper and lower value functions, and principle of optimality}

We now introduce the following notations: For $r>0$,
$$\cU_i[t,T;r]=\Big\{u_i\in\cU_i[t,T]\Bigm|\int_t^T|u_i(s)|^{\rho_i}ds\le
r\Big\},\qq i=1,2,$$
and
$$\left\{\ba{ll}
\ns\ds\cA_1[t,T;r]=\Big\{\a_1:\cU_2^1[t,T]\to\cU_1[t,T;r]\bigm|\a_1\in\cA_1[t,T]\Big\},\\
\ns\ds\cA_2[t,T;r]=\Big\{\a_2:\cU_1^1[t,T]\to\cU_2[t,T;r]\bigm|\a_2\in\cA_2[t,T]\Big\}.\ea\right.$$
We point out that although the upper and lower value functions are
formally defined in (\ref{upper-lower}), there seems to be no
guarantee that they are well-defined. The following result states
that under suitable conditions, $V^{\pm}(\cd\,,\cd)$ are indeed
well-defined.

\ms

\bf Theorem 5.4. \sl Let {\rm(H1)$'$--(H3)$'$ hold. Then the upper
and lower value functions $V^\pm(\cd\,,\cd)$ are well-defined and
there exists a constant $C>0$ such that
\bel{V}|V^\pm(t,x)|\le C\lan
x\ran{}^\m,\qq(t,x)\in[0,T]\times\dbR^n.\ee
Moreover,
\bel{5.10}\left\{\ba{ll}
\ns\ds
V^+(t,x)=\sup_{\a_2\in\cA_2[t,T;N(|x|)]}\inf_{u_1(\cd)\in\cU_1[t,T;N(|x|)]}
J(t,x;u_1(\cd),\a_2[u_1(\cd)]),\\
\ns\ds
V^-(t,x)=\inf_{\a_1\in\cA_1[t,T;N(|x|)]}\sup_{u_2(\cd)\in\cU_2[t,T;N(|x|)]}
J(t,x;\a_1[u_2(\cd)],u_2(\cd)),\ea\right.\ee
where $N(|x|)=C\lan x\ran{}^\m$, for some constant $C>0$.

\ms

\it Proof. \rm First of all, for any $(t,x)\in[0,T]\times\dbR^n$ and
$u_1(\cd)\in\cU_1[t,T]$, by Proposition 5.2, we have
$$\ba{ll}
\ns\ds\lan y(s)\ran\le C_0\Big\{\lan
x\ran+\int_t^s|u_1(r)|^{\si_1}dr\Big\}\le C_0\Big\{\2n\lan
x\ran+\|u_1(\cd)\|_{L^{\si_1}(t,T)}^{\si_1}\Big\}.\ea$$
Then
$$\ba{ll}
\ns\ds J(t,x;u_1(\cd),0)=\int_t^Tg(s,y(s),u_1(s),0)ds+h(y(T))\\
\ns\ds\ge\int_t^T\[c|u_1(s)|^{\rho_1}-L\lan y(s)\ran{}^\m\]ds
-L\lan y(T)\ran{}^\m\\
\ns\ds\ge\2n\int_t^T\2n\[c|u_1(s)|^{\rho_1}\2n -\1n
LC_0^\m\Big(\1n\lan x\ran\1n+\2n\int_t^s\2n|u_1(r)|^{\si_1}dr
\Big)^\m\]ds-LC_0^\m\(\lan
x\ran\1n+\|u_1(\cd)\|_{L^{\si_1}(t,T)}^{\si_1}\)^\m\\
\ns\ds\ge-C\lan
x\ran{}^\m\2n-C\|u_1(\cd)\|_{L^{\si_1}(t,T)}^{\si_1\m}\2n
+\int_t^Tc|u_1(s)|^{\rho_1}ds.\ea$$
Since (note $\m\ge1$)
$$\|u_1(\cd)\|_{L^{\si_1}(t,T)}^{\si_1\m}
=\(\int_t^T|u_1(r)|^{\si_1}dr\)^\m\le(T-t)^{\m-1}
\int_t^T|u_1(r)|^{\si_1\m}dr,$$
we obtain (taking into account $\si_1\m<\rho_1$)
\bel{5.9}\ba{ll}
\ns\ds J(t,x;u_1(\cd),0)\ge-C\lan
x\ran{}^\m+\int_t^T\[c|u_1(s)|^{\rho_1}-C|u_1(s)|^{\si_1\m}\]ds\\
\ns\ds\qq\qq\qq\;\ge-C\lan
x\ran{}^\m+{c\over2}\int_t^T|u_1(s)|^{\rho_1}ds\ge-C\lan
x\ran{}^\m.\ea\ee
Consequently,
$$\ba{ll}
\ns\ds
V^+(t,x)=\sup_{\a_2\in\cA_2[t,T]}\inf_{u_1(\cd)\in\cU_1[t,T]}J(t,x;u_1(\cd),\a_2[u_1(\cd)])\\
\ns\ds\qq\qq\ge\inf_{u_1(\cd)\in\cU_1[t,T]}J(t,x;u_1(\cd),\a_2^0[u_1(\cd)])\ge-C\lan
x\ran{}^\m.\ea$$
Likewise, for any $u_2(\cd)\in\cU_2[t,T]$, we have
\bel{5.12}J(t,x;0,u_2(\cd))=\int_t^Tg(s,y(s),0,u_2(s))ds+h(y(T))\le
C\lan x\ran{}^\m.\ee
Thus,
$$\ba{ll}
\ns\ds V^+(t,x)=\sup_{\a_2\in\cA_2[0,T]}\inf_{u_1(\cd)\in\cU_1[t,T]}
J(t,x;u_1(\cd),\a_2[u_1(\cd)])\\
\ns\ds\qq\qq\le\sup_{\a_2\in\cA_2[t,T]}J(t,x;u^0_1(\cd),\a_2[u^0_1(\cd)])\le
C\lan x\ran{}^\m.\ea$$
Similar results also hold for the lower value function
$V^-(\cd\,,\cd)$. Therefore, we obtain that $V^{\pm}(t,x)$ are
well-defined for all $(t,x)\in[0,T]\times\dbR^n$ and (\ref{V})
holds.

\ms

Next, for the constant $C>0$ appearing in (\ref{V}), we set
$$N(r)={4C\over c}\lan r\ran{}^\m.$$
Then for any $u_1(\cd)\in\cU_1[t,T]\setminus\cU_1[t,T;N(|x|)]$, from
(\ref{5.9}), we see that
$$\ba{ll}
\ns\ds J(t,x;u_1(\cd),\a_2^0[u_1(\cd)])\ge-C\lan
x\ran{}^\m+{c\over2}\int_t^T|u_1(s)|^{\rho_1}ds>
C\lan x\ran{}^\m\\
\ns\ds\ge
V^+(t,x)=\sup_{\a_2\in\cA_2[t,T]}\inf_{u_1(\cd)\in\cU_1[t,T]}J(t,x;u_1(\cd),\a_2[u_1(\cd)]).
\ea$$
Thus,
\bel{}V^+(t,x)=\sup_{\a_2\in\cA_2[t,T]}\inf_{u_1(\cd)\in\cU_1[t,T;N(|x|)]}J(t,x;u_1(\cd),\a_2[u_1(\cd)]).\ee
Consequently, from (\ref{5.12}), for any
$u_1(\cd)\in\cU_1[t,T;N(|x|)]$, we have
$$\ba{ll}
\ns\ds-C\lan x\ran{}^\m\le
V^+(t,x)\le\sup_{\a_2\in\cA_2[t,T]}J(t,x;u_1(\cd),\a_2[u_1(\cd)])\\
\ns\ds\qq\qq\le
C\lan x\ran{}^\m+C\int_t^T|u_1(s)|^{\rho_1}ds-{c\over2}\int_t^T|\a_2[u_1(\cd)](s)|^{\rho_2}ds\\
\ns\ds\qq\qq\le C\lan x\ran{}^\m+2C^2\lan
x\ran{}^\m-{c\over2}\int_t^T|\a_2[u_1(\cd)](s)|^{\rho_2}ds.\ea$$
This implies that
\bel{}{c\over2}\int_t^T|\a_2[u_1(\cd)](s)|^{\rho_2}ds\le\wt C\lan
x\ran{}^\m,\qq\forall u_1(\cd)\in\cU_1[t,T;N(|x|)],\ee
with $\wt C=2C(C+1)>0$ being another absolute constant. Hence, if we
replace the original $N(r)$ by the following:
$$N(r)={4\wt C\over c}\lan r\ran{}^\m,$$
and let
$$\cA_2[t,T;r]=\Big\{\a_2\in\cA_2[t,T]\Bigm|\int_t^T|\a_2[u_1(\cd)](s)|^{\rho_2}ds\le
N(|x|)\Big\},$$
then the first relation in (\ref{5.10}) holds.

\ms

The second relation in (\ref{5.10}) can be proved similarly. \endpf

\ms

Next, we want to establish a modified Bellman's principle of
optimality. To this end, for any $(t,x)\in[0,T)\times\dbR^n$ and
$\bar t\in(t,T]$, let
$$\cU_i[t,\bar t;r]=\Big\{u_i(\cd)\in\cU_i[t,T]\bigm|\int_t^{\bar t}|u_i(s)|^{\rho_i}ds\le r\Big\},\qq
i=1,2,$$
and
$$\left\{\ba{ll}
\ns\ds\cA_1[t,\bar t;r]=\Big\{\a_1:\cU_2^1[t,T]\to\cU_1[t,\bar t;r]\bigm|\a_1\in\cA_1[t,T]\Big\},\\
\ns\ds\cA_2[t,\bar t;r]=\Big\{\a_2:\cU_1^1[t,T]\to\cU_2[t,\bar
t;r]\bigm|\a_2\in\cA_2[t,T]\Big\}.\ea\right.$$
It is clear that
$$\left\{\ba{ll}
\ns\ds\cU_i[t,T;r]\subseteq\cU_i[t,\bar t;r]\subseteq\cU_i[t,T],\\
\ns\ds\cA_i[t,T;r]\subseteq\cA_i[t,\bar
t;r]\subseteq\cA_i[t,T],\ea\right.\qq i=1,2.$$
Thus, from the proof of Theorem 5.4, we see that for a suitable
choice of $N(\cd)$, say, $N(r)=C(1+r^\m)$ for some large $C>0$, the
following holds:
\bel{5.15}\left\{\ba{ll}
\ns\ds
V^+(t,x)=\sup_{\a_2\in\cA_2[t,\bar t;N(|x|)]}\inf_{u_1(\cd)\in\cU_1[t,T;N(|x|)]}J(t,x;u_1(\cd),\a_2[u_1(\cd)]),\\
\ns\ds V^-(t,x)=\inf_{\a_1\in\cA_1[t,\bar
t;N(|x|)]}\sup_{u_2(\cd)\in\cU_2[t,\bar
t;N(|x|)]}J(t,x;\a_1[u_2(\cd)],u_2(\cd)).\ea \right.\ee
We now state the following modified Bellman's principle of
optimality.

\ms

\bf Theorem 5.5. \sl Let {\rm(H1)$'$--(H3)$'$} hold. Let
$(t,x)\in[0,T)\times\dbR^n$ and $\bar t\in(t,T]$. Let
$N:[0,\infty)\to[0,\infty)$ be a nondecreasing continuous function
such that $(\ref{5.15})$ holds. Then
\bel{5.21}\ba{ll}
\ns\ds V^+(t,x)=\sup_{\a_2\in\cA_2[t,\bar
t;N(|x|)]}\inf_{u_1(\cd)\in\cU_1[t,\bar t;N(|x|)]}\Big\{\int_t^{\bar
t}g(s,y(s),u_1(s),\a_2[u_1(\cd)](s))ds+V^+(\bar t,y(\bar
t))\Big\},\ea\ee
and
\bel{5.22}\ba{ll}
\ns\ds V^-(t,x)=\inf_{\a_1\in\cA_1[t,\bar
t;N(|x|)]}\sup_{u_2(\cd)\in\cU_2[t,\bar t;N(|x|)]}\Big\{\int_t^{\bar
t}g(s,y(s),\a_1[u_2(\cd)](s),u_2(s))ds+V^-(\bar t,y(\bar
t))\Big\}.\ea\ee

\ms

\rm

We note that if in (\ref{5.21}) and (\ref{5.22}), $\cA_i[t,\bar
t;N(|x|)]$ and $\cU_i[t,\bar t;N(|x|)]$ are replaced by $\cA_i[t,T]$
and $\cU_i[t,T]$, respectively, the result is standard and the proof
is routine. However, in the above case, some careful modification is
necessary. For readers' convenience, we provide a proof in the
appendix.

\ms

We point out that our modified principle of optimality will play an
essential role in the next subsection.

\subsection{Continuity of upper and lower value functions}

In this subsection, we are going to establish the continuity of the
upper and lower value functions. Let us state the main results
now.

\ms

\bf Theorem 5.6. \sl Let {\rm(H1)$'$--(H3)$'$} hold. Then
$V^\pm(\cd\,,\cd)$ are continuous. Moreover, there exists a constant
$C>0$ and a nondecreasing continuous function
$N:[0,\infty)\to[0,\infty)$ such that the following estimates hold:
\bel{V-V}\ba{ll}
\ns\ds|V^\pm(t,x)-V^\pm(t,\bar x)|\le C\big(\lan x\ran\vee\lan
x\ran\big)^{\m-1}|x-\bar x|,\qq t\in[0,T],~x,\bar x\in\dbR^n,\ea\ee
and
\bel{5.19}|V^\pm(t,x)-V^\pm(\bar t,x)|\le N(|x|)|t-\bar
t|^{{\rho_1-\si_1\over\rho_1}\land{\rho_2-\si_2\over\rho_2}},\qq\forall
t,\bar t\in[0,T],~x\in\dbR^n.\ee

\it Proof. \rm We will only prove the conclusions for
$V^+(\cd\,,\cd)$. The conclusions for $V^-(\cd\,,\cd)$ can be proved
similarly.

\ms

First, let $0\le t\le T$, $x,\bar x\in\dbR^n$, and let $N(r)=C\lan
r\ran{}^\m$ for some $C>0$, such that (\ref{5.10}) holds. Take
\bel{5.14}u_1(\cd)\in\cU^{\rho_1}_1[t,T;N(|x|\vee|\bar
x|)],\q\a_2\in\wt\cA_2^{\rho_2}[t,T;N(|x|\vee|\bar x|)].\ee
Denote $u_2(\cd)=\a_2[u_1(\cd)]$. Then
$$\int_t^T|u_i(r)|^{\si_i}dr\le C\(\int_t^T|u_i(r)|^{\rho_i}dr\)^{\si_i\over\rho_i}
\le C\big(\1n\lan x\ran\vee\lan x\ran\big)^{\si_i\m\over\rho_i}\le
C\lan x\ran\vee\lan\bar x\ran,\qq i=1,2.$$
Making use of Proposition 5.1, we have
$$\ba{ll}
\ns\ds|y_{t,x}(s)|,|y_{t,\bar x}(s)|\le C_0\[\lan
x\ran{}\vee\lan\bar
x\ran+\int_t^T\(|u_1(r)|^{\si_1}+|u_2(r)|^{\si_2}\)dr\]\le
C\big(\lan x\ran\vee\lan\bar x\ran\big),\qq s\in[t,T],\ea$$
and
$$\ba{ll}
\ns\ds|y_{t,x}(s)-y_{t,\bar x}(s)|\le C_0|x-\bar x|,\qq
s\in[t,T].\ea$$
Consequently,
$$\ba{ll}
\ns\ds|J(t,x;u_1(\cd),u_2(\cd))-J(t,\bar x;u_1(\cd),u_2(\cd))|\\
\ns\ds\le\int_t^T|g(s,y_{t,x}(s),u_1(s),u_2(s))-g(s,y_{t,\bar
x}(s),u(s))|ds+|h(y_{t,x}(T))-h(y_{t,\bar x}(T))|\\
\ns\ds\le\int_t^T L\(\big(\1n\lan y_{t,x}(s)\ran\vee\lan y_{t,\bar
x}(s)\ran\big)^{\m-1}+|u_1(s)|^{\rho_1(\m-1)\over\m}+|u_2(s)|^{\rho_2(\m-1)\over\m}\)
|y_{t,x}(s)-y_{t,\bar x}(s)|ds\\
\ns\ds\q+L\(\lan y_{t,x}(T)\ran\vee\lan y_{t,\bar x}(T)\ran\)^{\m-1}
|y_{t,x}(T)-y_{t,\bar x}(T)|\\
\ns\ds\le C\Big\{\big(\lan x\ran\vee\lan\bar
x\ran\big)^{\m-1}+\(\int_t^T|u_1(s)|^{\rho_1}ds\)^{\m-1\over\m}
+\(\int_t^T|u_2(s)|^{\rho_2}ds\)^{\m-1\over\m}\Big\}|x-\bar x|\\
\ns\ds\le C\big(\lan x\ran\vee\lan\bar x\ran\big)^{\m-1}|x-\bar
x|.\ea$$
Since the above estimate is uniform in $(u_1(\cd),\a_2)$ satisfying
(\ref{5.14}), we obtain (\ref{V-V}) for $V^+(\cd\,,\cd)$.

\ms

We now prove the continuity in $t$. From the modified principle of
optimality, we see that for any $\e>0$, there exists an
$\a_2^\e\in\cA_2[t,\bar t;N(|x|)]$ such that
$$\ba{ll}
\ns\ds V^+(t,x)-\e\le\inf_{u_1(\cd)\in\cU_1[t,\bar
t;N(|x|)]}\Big\{\int_t^{\bar
t}g(s,y(s),u_1(\cd),\a_2^\e[u_1(\cd)](s))ds+V^+(\bar t,y(\bar t))\Big\}\\
\ns\ds\le\int_t^{\bar
t}g(s,y(s),0,\a_2^\e[u_1^0(\cd)](s))ds+V^+(\bar
t,y(\bar t))\\
\ns\ds\le\int_t^{\bar t}L\(\lan
y(s)\ran{}^\m-c|\a_2[u_1^0(\cd)](s))|^{\rho_2}\)ds+V^+(\bar
t,x)+|V^+(\bar
t,y(\bar t))-V^+(\bar t,x)|\\
\ns\ds\le\int_t^{\bar t}L\lan y(s)\ran{}^\m ds+V^+(\bar
t,x)+|V^+(\bar t,y(\bar t))-V^+(\bar t,x)|.\ea$$
By Proposition 5.2, we have (denote
$u_2^\e(\cd)=\a_2^\e[u_1^0(\cd)]$)
$$\ba{ll}
\ns\ds|y(\bar t)-x|\le C\[\lan x\ran(\bar t-t)+\int_t^{\bar
t}|u_2^\e(s)|^{\si_2}ds\]\\
\ns\ds\le C\[\lan x\ran(\bar t-t)+\(\int_t^{\bar
t}|u_2^\e(s)|^{\rho_2}ds\)^{\si_2\over\rho_2}(\bar
t-t)^{\rho_2-\si_2\over\rho_2}\]\le C\[\lan x\ran(\bar
t-t)+N(|x|)(\bar t-t)^{\rho_2-\si_2\over\rho_2}\].\ea$$
Also,
$$|y(s)|\le C_0\[\lan x\ran+\int_t^{\bar t}|u_2^\e(s)|^{\si_2}ds\]\le
N(|x|),\qq s\in[t,\bar t].$$
Hence, by the proved (\ref{V-V}), we obtain
$$|V^+(\bar t,y(\bar t))-V^+(\bar t,x)|\le N(|x|\vee|y(\bar
t))|y(\bar t)-x|\le N(|x|)(\bar t-t)^{\rho_2-\si_2\over\rho_2}.$$
Consequently,
$$V^+(t,x)-V^+(\bar t,x)\le N(|x|)(\bar
t-t)^{\rho_2-\si_2\over\rho_2}+\e,$$
which yields
$$V^+(t,x)-V^+(\bar t,x)\le N(|x|)(\bar
t-t)^{\rho_2-\si_2\over\rho_2}.$$
On the other hand,
$$V^+(t,x)\ge\inf_{u_1(\cd)\in\cU_1[t,T;N(|x|)]}\Big\{\int_t^{\bar
t}g(s,y(s),u_1(s),0)ds+V^+(\bar t,y(\bar t))\Big\}.$$
Hence, for any $\e>0$, there exists a
$u_1^\e(\cd)\in\cU_1[t,T;N(|x|)]$ such that
$$\ba{ll}
\ns\ds V^+(t,x)+\e\ge\int_t^{\bar t}g(s,y(s),u_1^\e(s),0)ds+V^+(\bar
t,y(\bar t))\\
\ns\ds\ge-\int_t^{\bar t}L\lan y(s)\ran{}^\m ds+c\int_t^{\bar
t}|u_1^\e(s)|^{\rho_1}ds+V^+(\bar t,x)-|V^+(\bar t,y(\bar
t))-V^+(\bar t,x)|\\
\ns\ds\ge-\int_t^{\bar t}L\lan y(s)\ran{}^\m ds+V^+(\bar
t,x)-|V^+(\bar t,y(\bar t))-V^+(\bar t,x)|.\ea$$
Now, in the current case, we have
$$\ba{ll}
\ns\ds|y(\bar t)-x|\le C\[\lan x\ran(\bar t-t)+\int_t^{\bar
t}|u_1^\e(s)|^{\si_1}ds\]\\
\ns\ds\le C\[\lan x\ran(\bar t-t)+\(\int_t^{\bar
t}|u_1^\e(s)|^{\rho_1}ds\)^{\si_1\over\rho_1}(\bar
t-t)^{\rho_1-\si_1\over\rho_1}\]\le C\[\lan x\ran(\bar
t-t)+N(|x|)(\bar t-t)^{\rho_1-\si_1\over\rho_1}\].\ea$$
Also,
$$|y(s)|\le C_0\[\lan x\ran+\int_t^{\bar t}|u_1^\e(s)|^{\si_1}ds\]\le
N(|x|),\qq s\in[t,\bar t].$$
Hence, by the proved (\ref{V-V}), we obtain
$$|V^+(\bar t,y(\bar t))-V^+(\bar t,x)|\le N(|x|\vee|y(\bar
t))|y(\bar t)-x|\le N(|x|)(\bar t-t)^{\rho_1-\si_1\over\rho_1}.$$
Consequently,
$$V^+(t,x)-V^+(\bar t,x)\ge-N(|x|)(\bar
t-t)^{\rho_1-\si_1\over\rho_1}-\e,$$
which yields
$$V^+(t,x)-V^+(\bar t,x)\ge-N(|x|)(\bar
t-t)^{\rho_1-\si_1\over\rho_1}.$$
Hence, we obtain the estimate (\ref{5.19}) for $V^+(\cd\,,\cd)$.
\endpf

\ms

\subsection{Characterization of the upper and lower value functions}

Having the above preparations, we are now at the position to
characterize the upper and the lower value functions of our
differential game. Recall that in order Theorem 4.3 applies, we need
the conditions (\ref{ml})--(\ref{h<}) (for the maps
$H(\cd\,,\cd\,,\cd)$ and $h(\cd)$ stated in (HJ) hold, and the upper
and lower value functions have to be Lipschitz continuous in a
particular form (see (\ref{|V-V|})). It is clear that the only thing
that we need is the compatibility condition (\ref{ml}) for the
numbers $\l_i,\n_i$ appeared in (\ref{3.24}) with the parameter $\m$
appeared in (H2) and (H2)$'$. Let us now look at what we need here.
From (\ref{3.24}) (which is for the upper value function
$V^+(\cd\,,\cd)$ only), and the similar set of conditions for lower
value function $V^-(\cd\,,\cd)$, we should require:
\bel{le1}\left\{\ba{ll}
\ns\ds{\si_1\m\over\rho_1}\le1,\q{\si_2\m\over\rho_2}\le1,
\q{(\m-1)\si_1\over\rho_1-\si_1}\le1,\q{(\m-1)\si_2\over\rho_2-\si_2}\le1,\\
\ns\ds{(\m-1)\si_1\rho_2\over
\rho_1(\rho_2-\si_2)}\le1,\q{(\m-1)\si_2\rho_1\over\rho_2(\rho_1-\si_1)}\le1,\\
\ns\ds{\si_1\si_2\m\over\rho_1\rho_2}+{(\m-1)\si_1\over\rho_1}\le1,\qq
{\si_1\si_2\m\over\rho_1\rho_2}+{(\m-1)\si_2\over\rho_2}\le1,\\
\ns\ds{\si_1\si_2\over
\rho_1\rho_2}+{(\m-1)\si_2(\si_1+\rho_1)\over\rho_1\rho_2}\le1,\qq
{\si_1\si_2\over
\rho_1\rho_2}+{(\m-1)\si_1(\si_2+\rho_2)\over\rho_1\rho_2}\le1,\\
\ns\ds{(\m-1)\si_1\si_2\over\rho_1(\rho_2-\si_2)}+{(\m-1)\si_2\over\rho_2}\le1,\qq
{(\m-1)\si_1\si_2\over\rho_2(\rho_1-\si_1)}+{(\m-1)\si_1\over\rho_1}\le1.
\ea\right.\ee
We now have the following proposition.

\ms

\bf Proposition 5.7. \sl Let
\bel{sm<r}\m\si_i\le\rho_i,\qq i=1,2.\ee
Then all the inequalities in $(\ref{le1})$ hold.

\ms

\it Proof. \rm First of all, we have that
$${\m\si_i\over\rho_i}\le1\q\iff\q{(\m-1)\si_i\over\rho_i-\si_i}\le1.$$
Thus, under (\ref{sm<r}), the last two inequalities in the first
line of (\ref{le1}) hold. Next, by the above equivalence and
$\m\ge1$,
$${(\m-1)\si_1\rho_2\over\rho_1(\rho_2-\si_2)}\le{(\m-1)\rho_2\over
\m(\rho_2-\si_2)}\le1,$$
and
$${(\m-1)\si_2\rho_1\over\rho_2(\rho_1-\si_1)}\le{(\m-1)\rho_1\over
\m(\rho_1-\si_1)}\le1.$$
Thus, the inequalities in the second line of (\ref{le1}) hold. Now,
for the third line, we have
$${\si_1\si_2\m\over\rho_1\rho_2}+{(\m-1)\si_1\over\rho_1}\le{\si_1\over\rho_1}
+{(\m-1)\si_1\over\rho_1}={\m\si_1\over\rho_1}\le1,$$
and
$${\si_1\si_2\m\over\rho_1\rho_2}+{(\m-1)\si_2\over\rho_2}\le{\si_2\over\rho_2}
+{(\m-1)\si_2\over\rho_2}={\m\si_2\over\rho_2}\le1.$$
This shows that the inequalities in the third line of (\ref{le1})
hold. We now look at the fourth line. It is seen that
$${\si_1\si_2\over
\rho_1\rho_2}+{(\m-1)\si_2(\si_1+\rho_1)\over\rho_1\rho_2}\le{\si_1\si_2\over
\rho_1\rho_2}+{(\m-1)\si_2(\si_1+\m\si_1)\over\rho_1\rho_2}
={\m^2\si_1\si_2\over\rho_1\rho_2}\le1,$$
and
$${\si_1\si_2\over
\rho_1\rho_2}+{(\m-1)\si_1(\si_2+\rho_2)\over\rho_1\rho_2}\le{\si_1\si_2\over
\rho_1\rho_2}+{(\m-1)\si_1(\si_2+\m\si_2)\over\rho_1\rho_2}={\m^2\si_1\si_2
\over\rho_1\rho_2}\le1.$$
Finally, for the fifth line, we have (making use of the inequalities
in the second line of (\ref{le1}))
$${(\m-1)\si_1\si_2\over\rho_1(\rho_2-\si_2)}+{(\m-1)\si_2\over\rho_2}
={\si_2\over\rho_2}\[{(\m-1)\si_1\rho_2\over\rho_1(\rho_2-\si_2)}+\m-1\]
\le{\si_2\m\over\rho_2}\le1,$$
and
$${(\m-1)\si_1\si_2\over\rho_2(\rho_1-\si_1)}+{(\m-1)\si_1\over\rho_1}
={\si_1\over\rho_1}\[{(\m-1)\si_2\rho_1\over\rho_2(\rho_1-\si_1)}+\m-1\]
\le{\si_1\m\over\rho_1}\le1.$$
This completes the proof. \endpf

\ms

With the above result, we have the following theorem.

\ms

\bf Theorem 5.8. \sl Let {\rm(H1)$'$--(H3)$'$} hold. Then
$V^\pm(\cd\,,\cd)$ are the unique viscosity solution to the upper
and lower HJI equations $(\ref{HJI0})$, respectively. Further, if
the Isaacs' condition holds:
\bel{}H^+(t,x,p)=H^-(t,x,p),\qq\forall(t,x,p)\in[0,T]\times\dbR^n\times\dbR^n,\ee
then
\bel{V=V}V^+(t,x)=V^-(t,x),\qq\forall(t,x)\in[0,T]\times\dbR^n.\ee

\ms

\section{Remarks on the Existence of Viscosity Solutions to HJ
Equations.}

\rm

We have seen that under (H1)--(H3), the upper and lower Hamiltonians
can be well-defined and the corresponding upper and lower HJI
equations can be well-formulated. Moreover, we have proved the
uniqueness of the viscosity solutions to the upper and lower HJI
equations within a suitable class of locally Lipschitz continuous
functions. On the other hand, we have introduced a little stronger
hypotheses (H1)$'$--(H3)$'$ to obtain the upper and lower value
functions $V^\pm(\cd\,,\cd)$ being well-defined so that the
corresponding upper and lower HJI equations have viscosity
solutions. In another word, weaker conditions ensure the uniqueness
of viscosity solutions to the upper and lower HJI equations, and
stronger conditions seem to be needed for the existence. There are
some general existence results of viscosity solutions for the first
order HJ equations in the literature, see \cite{Lions 1982, Barles
1984, Souganidis 1985, Friedman-Soudanidis 1986, Crandall-Lions
1987}. A natural question is whether the conditions that we assumed
for the existence of viscosity solutions are sharp (or close to be
necessary). In this section, we present a simple situation which
tells us that our conditions are sharp in some sense.

\ms

We consider the following one-dimensional controlled linear system:
\bel{}\left\{\ba{ll}
\ns\ds\dot y(s)=Ay(s)+B_1u_1(s)+B_2u_2(s),\qq s\in[t,T],\\
\ns\ds y(t)=x,\ea\right.\ee
with the performance functional:
\bel{}J(t,x;u_1(\cd),u_2(\cd))=\int_t^T\[Qy(s)^2+R_1u_1(s)^2-R_2u_2(s)^2\]ds+Gy(T)^2,\ee
where $A,B_1,B_2,A,R_1,R_2,G\in\dbR$. We assume that
\bel{6.3}R_1,R_2>0.\ee
Note that in the current case,
$$\si_1=\si_2=1,\q\m=\rho_1=\rho_2=2.$$
Thus,
$$\m\si_i=\rho_i,\qq i=1,2,$$
which violates (\ref{5.4}). In the current case, we have
\bel{}\ba{ll}
\ns\ds H^\pm(t,x,p)=H(t,x,p)=\inf_{u_1}\sup_{u_2}\[pf(t,x,u_1,u_2)+g(t,x,u_1,u_2)\]\\
\ns\ds=Apx+Qx^2+\inf_{u_1}\[R_1u_1^2+pB_1u_1\]-\inf_{u_2}\[R_2u_2^2-pB_2u_2\]
=Apx+Qx^2+\({B_2^2\over4R_2}-{B_1^2\over4R_1}\)p^2.\ea\ee
Consequently, the upper and lower HJI equation have the same form:
\bel{6.5}\left\{\ba{ll}
\ns\ds
V_t(t,x)+AxV_x(t,x)+Qx^2+\({B_2^2\over4R_2}-{B_1^2\over4R_1}\)V_x(t,x)^2=0,\q(t,x)\in[0,T]\times\dbR,\\
\ns\ds V(T,x)=Gx^2,\qq x\in\dbR.\ea\right.\ee
If the above HJI equation has a viscosity solution, by the
uniqueness, the solution has to be of the following form:
\bel{6.6}V(t,x)=p(t)x^2,\qq(t,x)\in[0,T]\times\dbR,\ee
where $p(\cd)$ is the solution to the following Riccati equation:
\bel{Riccati-1}\left\{\ba{ll}
\ns\ds\dot
p(t)+2Ap(t)+Q+\({B_2^2\over R_2}-{B_1^2\over R_1}\)p(t)^2=0,\qq t\in[0,T],\\
\ns\ds p(T)=G.\ea\right.\ee
In another word, the solvability of (\ref{6.5}) is equivalent to
that of (\ref{Riccati-1}).

\ms

Our claim is that Riccati equation (\ref{Riccati-1}) is not always
solvable for any $T>0$. To state our result in a relatively neat
way, let us rewrite equation (\ref{Riccati-1}) as follows:
\bel{6.8}\left\{\ba{ll}
\ns\ds\dot p+\a p+\b p^2+\g=0,\\
\ns\ds p(T)=g,\ea\right.\ee
with
$$\a=2A,\qq\b={B_2^2\over R_2^2}-{B_1^2\over R_1^2},\qq\g=Q,\qq g=G.$$
Note that $\b$ could be positive, negative, or zero. We have the
following result.

\ms

\bf Proposition 6.1. \sl Riccati equation $(\ref{6.8})$ admits a
solution on $[0,T]$ for any $T>0$ if and only if one of the
following holds:
\bel{6.9}\a^2-4\b\g\ge0,\q 2\b g+\a-\sqrt{\a^2-4\b\g}\le0;\ee

\rm

The proof is elementary and straightforward. For reader's
convenience, we provide a proof in the appendix.

\ms

It is clear that there are a lot of cases for which the Riccati
equation is not solvable. For example,
$$\a=\b=\g=1,$$
which violates (\ref{6.9}). Also, the case
$$\a=0,\q\b=-1,\q\g=1,\q g=-2,$$
which also violates (\ref{6.9}). For the above two cases, Riccati
equation (\ref{6.8}) does not have a global solution on $[0,T]$ for
some $T>0$. Correspondingly we have some two-person zero-sum
differential game with unbounded controls for which the coercivity
condition (\ref{5.4}) fails and the upper and lower value functions
could not be defined on the whole time interval $[0,T]$, or
equivalently, the corresponding upper/lower HJI equation have no
viscosity solutions on $[0,T]$.

\ms

\ms

\no\bf Acknowledgement. \rm The authors would like to thank the
referee for informing the authors several important references in
the field, especially some most recent papers. Also, some comments
made by Professor Y.~Hu (of University of Rennes 1, France) on the
previous version are really appreciated.

\rm

\bs

\no{\Large\bf Appendix}

\rm

\ms

\it Proof of Theorem 5.4. \rm We only prove (\ref{5.21}). The other
can be proved similarly. Since $N(|x|)$ and $\bar t$ are fixed, for
notational simplicity, we denote below that
$$\wt\cU_1=\cU_1[t,\bar t;N(|x|)],\qq\wt\cA_2=\cA_2[t,\bar t;N(|x|)].$$
Denote the right hand side of (\ref{5.21}) by $\h V^+(t,x)$. For any
$\e>0$, there exists an $\a_2^\e\in\wt\cA_2$ such that
$$\h V^+(t,x)-\e<\inf_{u_1(\cd)\in\wt\cU_1}\Big\{\int_t^{\bar
t}g(s,y(s),u_1(s),\a^\e_2[u_1(\cd)](s))ds+V^+(\bar t,y(\bar
t))\Big\}.$$
By the definition of $V^+(\bar t,y (\bar t))$, there exists an
$\bar\a_2^\e\in\cA_2[\bar t,T]$ such that
$$V^+(\bar t,y(\bar t))-\e<\inf_{\bar u_1(\cd)\in\cU_1[\bar
t,T]}J(\bar t,y(\bar t);\bar u_1(\cd),\bar\a_2^\e[\bar u_1(\cd)]).$$
Now, we define an extension $\h\a^\e_2\in\cA_2[t,T]$ of
$\a_2^\e\in\cA_2[\bar t,T]$ as follows: For any
$u_1(\cd)\in\cU_1[t,T]$,
$$\h\a^\e_2[u_1(\cd)](s)=\left\{\ba{ll}
\ns\ds\a_2^\e[u_1(\cd)](s),\qq\qq\qq s\in[t,\bar t),\\
\ns\ds\bar\a^\e_2[u_1(\cd)\big|_{[\bar t,T]}](s),\qq\qq s\in[\bar
t,T].\ea\right.$$
Since $\a_2^\e\in\wt\cA_2$, we have
$$\int_t^{\bar t}|\h\a^\e[u_1(\cd)](s)|^{\rho_2}ds=\int_t^{\bar t}|\a_2^\e[u_1(\cd)](s)|^{\rho_2}ds\le N(|x|).$$
This means that $\h\a_2^\e\in\wt\cA_2$. Consequently,
$$\ba{ll}
\ns\ds
V^+(t,x)\ge\inf_{u_1(\cd)\in\wt\cU_1}J(t,x;u_1(\cd),\h\a_2^\e[u_1(\cd)])\\
\ns\ds=\inf_{u_1(\cd)\in\wt\cU_1}\Big\{\int_t^{\bar
t}g(s,y(s),u_1(s),\a_2^\e[u_1(\cd)](s))ds+J(\bar t,y(\bar
t);u_1(\cd)\big|_{[\bar
t,T]},\bar\a_2^\e[u_1(\cd)\big|_{[\bar t,T]})\Big\}\\
\ns\ds\ge\inf_{u_1(\cd)\in\wt\cU_1}\Big\{\int_t^{\bar
t}g(s,y(s),u_1(s),\a_2^\e[u_1(\cd)](s))ds+\inf_{\bar
u_1(\cd)\in\cU_1[\bar t,T]}J(\bar t,y(\bar t);\bar
u_1(\cd),\bar\a_2^\e[\bar u_1(\cd))\Big\}\\
\ns\ds\ge\inf_{u_1(\cd)\in\wt\cU_1}\Big\{\int_t^{\bar
t}g(s,y(s),u_1(s),\a_2^\e[u_1(\cd)](s))ds+V^+(\bar t,y(\bar
t))\Big\}-\e\ge\h V^+(t,x)-2\e.\ea$$
Since $\e>0$ is arbitrary, we obtain
$$\h V^+(t,x)\le V^+(t,x).$$
On the other hand, for any $\e>0$, there exists an
$\a_2^\e\in\wt\cA_2$ such that
$$V^+(t,x)-\e<\inf_{u_1(\cd)\in\wt\cU_1}J(t,x;u_1(\cd),\a_2^\e[u_1(\cd)]).$$
Also, by definition of $\h V^+(t,x)$,
$$\h V^+(t,x)\ge\inf_{u_1(\cd)\in\wt\cU_1}\Big\{\int_t^{\bar
t}g(s,y(s),u_1(s),\a_2^\e[u_1(\cd)](s))ds+V^+(\bar t,y(\bar
t))\Big\}.$$
Thus, there exists a $u_1^\e(\cd)\in\wt\cU_1$ such that
$$\h V^+(t,x)+\e\ge\int_t^{\bar
t}g(s,y(s),u_1^\e(s),\a_2^\e[u_1^\e(\cd)](s))ds+V^+(\bar t,y(\bar
t)).$$
Now, for any $\bar u_1(\cd)\in\cU_1[\bar t,T]$, define a particular
extension $\wt u_1(\cd)\in\cU_1[t,T]$ by the following:
$$\wt u_1(s)=\left\{\ba{ll}
\ns\ds u_1^\e(s),\qq s\in[t,\bar t),\\
\ns\ds\bar u_1(s),\qq s\in[\bar t,T].\ea\right.$$
Namely, we patch $u_1^\e(\cd)$ to $\bar u_1(\cd)$ on $[t,\bar t)$.
Since
$$\int_t^{\bar t}|\wt u_1(s)|^{\rho_1}ds=\int_t^{\bar t}|u_1^\e(s)|^{\rho_1}ds\le N(|x|),$$
we see that $\wt u_1(\cd)\in\wt\cU_1$. Next, we define a restriction
$\bar\a_2^\e\in\cA[\bar t,T]$ of $\a_2^\e\in\wt\cA_2$, as follows:
$$\bar\a_2^\e[\bar u_1(\cd)]=\a_2^\e[\wt u_1(\cd)].$$
For such an $\bar\a_2^\e$, we have
$$V^+(\bar t,y(\bar t))\ge\inf_{\bar u_1(\cd)\in\cU_1[\bar
t,T]}J(\bar t,y(\bar t),\bar u_1(\cd),\bar\a_2^\e[\bar u_1(\cd)]).$$
Hence, there exists a $\bar u_1^\e(\cd)\in\cU_1[\bar t,T]$ such that
$$V^+(\bar t,y(\bar t))+\e>J(\bar t,y(\bar t),\bar
u_1^\e(\cd),\bar\a_2^\e[\bar u_1^\e(\cd)]).$$
Then we further let
$$\wt u_1^\e(s)=\left\{\ba{ll}
\ns\ds u_1^\e(s),\qq\q s\in[t,\bar t),\\
\ns\ds\bar u_1^\e(s),\qq\q s\in[\bar t,T].\ea\right.$$
Again, $\wt u_1^\e(\cd)\in\wt\cU_1$, and therefore,
$$\ba{ll}
\ns\ds\h V^+(t,x)+\e\ge\int_t^{\bar
t}g(s,y(s),u_1^\e(s),\a_2^\e[u_1^\e(\cd)](s))ds+V^+(\bar t,y(\bar
t))\\
\ns\ds\ge\int_t^{\bar
t}g(s,y(s),u_1^\e(s),\a_2^\e[u_1^\e(\cd)](s))ds+J(\bar t,y(\bar
t),\bar u_1^\e(\cd),\bar\a_2^\e[\bar u_1^\e(\cd)])-\e\\
\ns\ds=J(t,x;\wt u_1^\e(\cd),\a_2^\e[\wt
u_1^\e(\cd)])-\e\ge\inf_{u_1(\cd)\in\wt\cU_1[t,T]}J(t,x;u_1(\cd),\a_2^\e[u_1(\cd)])-\e\ge
V^+(t,x)-2\e.\ea$$
Since $\e>0$ is arbitrary, we obtain
$$\h V^+(t,x)\ge V^+(t,x).$$
This completes the proof. \endpf

\ms

\it Proof of Proposition 6.1. \rm Recall that we are considering the
following Riccati equation:
$$\left\{\ba{ll}
\ns\ds\dot p+\a p+\b p^2+\g=0,\\
\ns\ds p(T)=g,\ea\right.$$

\ms

\it Case 1. \rm $\b=0$. The Riccati equation reads
$$\left\{\ba{ll}
\ns\ds\dot p+\a p+\g=0,\\
\ns\ds p(T)=g.\ea\right.$$
This is an initial value problem for a linear equation, which admits
a unique global solution $p(\cd)$ on $[0,T]$.

\ms

\it Case 2. \rm $\b\ne0$. Then Riccati equation reads
$$\left\{\ba{ll}
\ns\ds\dot p+\b\[\(p+{\a\over2\b}\)^2+{4\b\g-\a^2\over4\b^2}\]=0,\\
\ns\ds p(T)=g.\ea\right.$$
Let
$$\k={\sqrt{|\a^2-4\b\g|}\over2|\b|}\ge0.$$
There are three subcases.

\ms

\it Subscase 1. \rm $\a^2-4\b\g=0$. The Riccati equation becomes
$$\left\{\ba{ll}
\ns\ds\dot p+\b\(p+{\a\over2\b}\)^2=0,\\
\ns\ds p(T)=g.\ea\right.$$
Therefore, in the case
$$2\b g+\a=0,$$
we have that $p(t)\equiv-{\a\over2\b}$ is the (unique) global
solution on $[0,T]$. Now, let
$$2\b g+\a\ne0.$$
Then we have
$${dp\over(p+{\a\over2\b})^2}=-\b dt,$$
which leads to
$${1\over p(t)+{\a\over2\b}}={1\over g+{\a\over2\b}}-\b(T-t)={2\b-\b(2\b g+\a)(T-t)\over2\b g+\a}.$$
Thus,
$$p(t)=-{\a\over2\b}+{2\b g+\a\over2\b-\b(2\b g+\a)(T-t)},$$
which is well-defined on $[0,T]$ if and only if
$$2-(2\b g+\a)(T-t)\ne0,\qq t\in[0,T].$$
This is equivalent to the following:
$$(2\b g+\a)T<2.$$
The above is true for all $T>0$ if and only if
$$2\b g+\a\le0,$$

\ms

\it Subcase 2. \rm $\a^2-4\b\g<0$. The Riccati equation is
$$\dot p+\b\[\(p+{\a\over2\b}\)^2+\k^2\]=0.$$
Hence,
$${dp\over(p+{\a\over2\b})^2+\k^2}=-\b dt,$$
which results in
$${1\over\k}\tan^{-1}\[{1\over\k}\(p(t)+{\a\over2\b}\)\]=-\b t+C.$$
By the terminal condition,
$$C=\b T+{1\over\k}\tan^{-1}\[{1\over\k}\(g+{\a\over2\b}\)\]$$
Consequently,
$$\tan^{-1}\[{1\over\k}\(p(t)+{\a\over2\b}\)\]=\k\b(T-t)+\tan^{-1}\[{1\over\k}\(g+{\a\over2\b}\)\].$$
Then
$$p(t)={\a\over2\b}+\k\tan\Big\{\k\b(T-t)+\tan^{-1}\({2\b g+\a\over2\k\b}\)\Big\}.$$
The above is well-defined for $t\in[0,T]$ if and only if
$$-{\pi\over2}<\tan^{-1}{2\b g+\a\over2\k\b}+\k\b T<{\pi\over2},$$
which is true for all $T>0$ if and only if $\b=0$.

\ms

\it Subcase 3. \rm $\a^2-4\b\g>0$. The Riccati equation becomes
$$\dot p+\b\[\(p+{\a\over2\b}\)^2-\k^2\]=0.$$
If
$$(2\b g+\a-2\k\b)(2\b
g+\a+2\k\b)\equiv4\b^2\(g+{\a\over2\b}-\k\)\(g+{\a\over2\b}+\k\)=0,\rq{\rm(A1)}$$
then one of the following
$$p(t)\equiv-{\a\over2\b}\pm\k,\qq t\in[0,T],$$
is the unique global solution to the Riccati equation. We now let
$$(2\b g+\a-2\k\b)(2\b
g+\a+2\k\b)\equiv4\b^2\(g+{\a\over2\b}-\k\)\(g+{\a\over2\b}+\k\)\ne0.$$
Then
$${dp\over(p+{\a\over2\b})^2-\k^2}=-\b dt.$$
Hence,
$${1\over2\k}\ln\Big|{p(t)+{\a\over2\b}-\k\over
p(t)+{\a\over2\b}+\k}\Big|=-\b t+\wt C,$$
which implies
$${p(t)+{\a\over2\b}-\k\over
p(t)+{\a\over2\b}+\k}=Ce^{-2\k\b t},$$
with
$$C=e^{2\k\b T}{g+{\a\over2\b}-\k\over
g+{\a\over2\b}+\k}=e^{2\k\b T}{2\b g+\a-2\k\b\over2\b g+\a+2\k\b}.$$
Then
$${p(t)+{\a\over2\b}-\k\over
p(t)+{\a\over2\b}+\k}=e^{2\k\b (T-t)}{2\b g+\a-2\k\b\over2\b
g+\a+2\k\b}.$$
Consequently,
$$p(t)+{\a\over2\b}-\k=e^{2\k\b
(T-t)}{2\b g+\a-2\k\b\over2\b g+\a+2\k\b}\[p(t)+{\a\over2\b}+\k\].$$
Thus, $p(\cd)$ globally exists on $[0,T]$ if and only if
$$e^{2\k\b(T-t)}{2\b g+\a-2\k\b\over2\b
g+\a+2\k\b}-1\ne0,\qq\forall t\in[0,T],$$
which is equivalent to
$$\psi(t)\equiv e^{2\k\b(T-t)}(2\b g+\a-2\k\b)-(2\b g+\a+2\k\b)\ne0,\qq\forall t\in[0,T].$$
Since $\psi'(t)$ does not change sign on $[0,T]$, the above is
equivalent to the following:
$$0<\psi(0)\psi(T)=\[e^{2\k\b T}(2\b g+\a-2\k\b)-(2\b g+\a+2\k\b)\](-4\k\b),$$
which is equivalent to
$$\[e^{2\k\b T}(2\b g+\a-2\k\b)-(2\b g+\a+2\k\b)\]\b<0.$$
Note when (A1) holds, the above it true. In the case $\b>0$, the
above reads
$$e^{2\k\b T}(2\b g+\a-2\k\b)<2\b g+\a+2\k\b,$$
which is true for all $T>0$ if and only if
$$2\b g+\a-2\k\b\le0.\rq{\rm(A2)}$$
Finally, if $\b<0$, then
$$\ba{ll}
\ns\ds0<e^{2\k\b T}(2\b g+\a-2\k\b)-(2\b g+\a+2\k\b)\\
\ns\ds\q=e^{-2\k|\b|T}(-2|\b|g+\a+2\k|\b|)-(-2|\b|g+\a-2\k|\b|)\\
\ns\ds\q=e^{-2\k|\b|T}\[-\(2|\b|g-\a-2\k|\b|\)+e^{2\k|\b|T}\(2|\b|g-\a+2\k|\b|\)\],\ea$$
which is true for all $T>0$ if and only if
$$0\le2|\b|g-\a+2\k|\b|=-(2\b g+\a-2\k|\b|).$$
Thus,
$$2\b g+\a-2\k|\b|\le0.$$
which has the same form as (A2). This completes the proof.
\endpf

\ms


\begin{thebibliography}{9}

\nobreak

\ms

\rm

\bibitem{Bardi-Capuzzo-Dolcetta 1997} M.~Bardi and
I.~Capuzzo-Dolcetta, \sl Optimal Control and Viscosity Solutions of
Hamilton-Jacobi-Bellman Equations, \rm Birkh\"auser, Boston, 1997.

\bibitem{Bardi-Da Lio 1997} M.~Bardi and F.~Da Lio, \it On the Bellman equation for some
unbounded control problems, \sl NoDEA, \rm 4 (1997), 491--510.

\bibitem{Barles 1984} G.~Barles, \it Existence results for first
order Hamilton-Jacobi equations, \sl Annales de l'I.H.P., \rm 1
(1984), 325--340.

\bibitem{Biton 2001} S.~Biton, \it Nonlinear monotone semigroups and
viscosity solutions, \sl Ann. I. H. Poincar\'e-AN, \rm 18 (2001),
383--402.

\bibitem{Crandall-Lions 1983} M.~G.~Crandall and P.~L.~Lions, \it
Viscosity solutions of Hamilton-Jacobi equations, \sl Trans. AMS,
\rm 277 (1983), 1--42.

\bibitem{Crandall-Lions 1986} M.~G.~Crandall and P.~L.~Lions, \it On
existence and uniqueness of solutions of Hamilton-Jacobi equations,
\sl Nonlinear Anal., \rm 10 (1986), 353--370.

\bibitem{Crandall-Lions 1987} M.~G.~Crandall and P.~L.~Lions, \it
Remarks on the existence and uniqueness of unbounded viscosity
solutions of Hamilton-Jacobi equations, \sl Illinois J. Math., \rm
31 (1987), 665--688.

\bibitem{Da Lio 2000} F.~Da Lio, \it On the Bellman equation for
infinite horizon problems with unblounded cost functional, \sl Appl.
Math. Optim., \rm 41 (2000), 171--197.

\bibitem{Da Lio--Ley 2006}
F.~Da Lio and O.~Ley, \it Uniqueness results for second-order
Bellman-Isaacs equations under quadratic growth assumptions and
applications, \sl SIAM J. Control Optim., \rm 45 (2006), 74--106.

\bibitem{Da Lio--Ley 2011} F.~Da Lio and O.~Ley, \it Convex Hamilton-Jacobi equations
under superlinear growth conditions on data, \sl Appl. Math. Optim.,
\rm 63 (2011), 309--339.

\bibitem{Elliott-Kalton 1972} R.~J.~Elliott and N.~J.~Kalton, \it
The existence of value in differential games, \sl Memoirs of AMS,
\rm No. 126. Amer. Math. Soc., Providence, R.I., 1972.

\bibitem{Evans-Souganidis 1984} L.~C.~Evans and P.~E.~Souganidis,
\it Differential games and representation formulas for solutions of
Hamilton-Jacobi-Isaacs equations, \sl Indiana Univ. Math. J., \rm 5
(1984), 773--797.

\bibitem{Fleming-Souganidis 1989} W.~H.~Fleming and
P.~E.~Souganidis, \it On the existence of value functions of
two-players, zero-sum stochastic differential games, \sl Indiana
Univ. Math. J., \rm 38 (1989), 293--314.

\bibitem{Friedman-Soudanidis 1986} A.~Friedman and P.~E.~Souganidis,
\it Blow-up solutions of Hamilton-Jacobi equations, \sl Comm. PDEs,
\rm 11 (1986), 397--443.

\bibitem{Ishii 1984} H.~Ishii, \it Uniqueness of unbounded viscosity
solutions of Hamilton-Jacobi equations, \sl Indiana Univ. Math. J.,
\rm 33 (1984), 721--748.

\bibitem{Ishii 1988} H.~Ishii, \it Representation of solutions of Hamilton-Jacobi
equations, \sl Nonlinear Anal., \rm 12 (1988), 121-–146.

\bibitem{Lions 1982} P.~L.~Lions, \sl Generalized Solutions of
Hamilton-Jacobi equations, \rm Pitman, London, 1982.

\bibitem{Lions-Souganidis 1985} P.~L.~Lions and P.~E.~Souganidis,
\it Differential games, optimal conrol and directional derivatives
of viscosity solutions of Bellman's and Isaacs' equations, \sl SIAM
J. Control Optim., \rm 23 (1985), 566--583.

\bibitem{McEneaney 1998} W.~McEneaney, \it A uniqueness result for the Isaacs equation
corresponding to nonlinear $H_\infty$ control, \sl Math. Control
Signals Systems, \rm 11 (1998), 303--334.


\bibitem{Rampazzo 1998} F.~Rampazzo, \it Differential games with unbounded versus bounded
controls, \sl SIAM J. Control Optim., \rm 36 (1998), 814-–839.

\bibitem{Soravia 1999} P.~Soravia, \it Equivalence between nonlinear
$\cH_\infty$ control problems and existence of viscosity solutions
of Hamilton-Jacobi-Isaacs equations, \sl Appl. Math. Optim., \rm 39
(1999), 17--32.

\bibitem{Soravia 2004} M.~Garavello and P.~Soravia, \it Optimality
principles and uniqueness for Bellman equations of unbounded control
problems with discontinuous running cost, \sl NoDEA, \rm 11 (2004),
271--298.

\bibitem{Souganidis 1985} P.~E.~Souganidis, \it Existence of
viscosity solution of Hamilton-Jacobi equations, \sl J. Diff. Eqs.,
\rm 56 (1985), 345--390.

\bibitem{Yong 1994} J.~Yong, \it Zero-sum differential games
involving impusle controls, \sl Appl. Math. Optim., \rm 29 (1994),
243--261.

\bibitem{You 2002} Y.~You, \it Syntheses of differential games and
pseudo-Riccati equations, \sl Abstr. Appl. Anal., \rm 7 (2002),
61--83.




\end{thebibliography}
\end{document}